\documentclass[preprints,article,accept,moreauthors,pdftex]{Definitions/mdpi} 

\firstpage{1} 
\pubvolume{1}
\issuenum{1}
\articlenumber{0}
\pubyear{2021}
\copyrightyear{2020}
\datereceived{} 
\dateaccepted{}
\datepublished{} 
\hreflink{https://doi.org/} % If needed use \linebreak
    \pdfoutput=1

\usepackage[USenglish]{babel} % hyphenation

\usepackage[capitalize]{cleveref}

\usepackage{bm} % for bold math symbols (including Greek letters)
\usepackage{nicefrac} % `dfrac` for nicer inline fractions

\usepackage{mathtools} % to use `\clap` when supersetting stuff over operators

\usepackage[ruled]{algorithm2e}
\SetKwInput{Configuration}{Configuration}

\Title{Derivative-Free Multiobjective Trust Region Descent Method Using Radial Basis Function Surrogate Models}

\TitleCitation{Derivative-Free Multiobjective Trust Region Descent Method Using Radial Basis Function Surrogate Models}

\Author{
    Manuel Berkemeier $^{1}$\orcidA{} and 
    Sebastian Peitz $^{2}$\orcidB{}
}

\AuthorNames{Manuel Berkemeier and Sebastian Peitz}

\AuthorCitation{Berkemeier, M.; Peitz, S.}
\address{%
$^{1}$ \quad Paderborn University; manuelbb@math.upb.de\\
$^{2}$ \quad Paderborn University; speitz@math.upb.de}

\corres{Correspondence: manuelbb@math.upb.de}

\abstract{
    We present a flexible trust region descend algorithm for unconstrained and 
    convexly constrained multiobjective optimization problems.
    It is targeted at heterogeneous and expensive problems, i.e., 
    problems that have at least one objective function that is computationally 
    expensive.
    The method is derivative-free in the sense that neither need derivative information be available for the 
    expensive objectives nor are gradients approximated using repeated function evaluations as 
    is the case in finite-difference methods.
    Instead, a multiobjective trust region approach is used that works similarly to its well-known scalar pendants.
    Local surrogate models constructed from evaluation data of the true objective functions 
    are employed to compute possible descent directions.
    In contrast to existing multiobjective trust region algorithms, these surrogates are not polynomial but 
    carefully constructed radial basis function networks.
    This has the important advantage that the number of data points scales linearly with the parameter space dimension.
    The local models qualify as \emph{fully linear} and the corresponding general scalar framework is 
    adapted for problems with multiple objectives.
    Convergence to Pareto critical points is proven and numerical examples illustrate our findings.
}

\keyword{
multiobjective optimization;
trust region methods;
multiobjective descent;
derivative-free optimization;
radial basis functions;
fully linear models
}

\Crefname{algocf}{Algorithm}{Algorithms}
\Crefname{assumption}{Assumption}{Assumptions}
\Crefname{Lemma}{Lemma}{Lemmata}
\Crefname{max_prop}{the maximum property}{the maximum properties}
\Crefname{Appendix}{Appendix}{Appendices}
\Crefname{specialtable}{Table}{Tables}

\newcommand{\expensiveindices}{I_{\text{ex}}}
\newcommand{\cheapinidces}{I_{\text{cheap}}}

\newcommand{\ve}[1]{%
    \ifcat\noexpand#1\relax % check if the argument is a control sequence
    \bm{#1}% probably Greek
    \else
    \bm{\mathrm{#1}}% single character
    \fi%
}
\newcommand{\ita}[1]{^{(#1)}}
\newcommand{\itat}{\ita{t}}

\DeclareMathOperator*{\argmin}{arg\,min}

\newcommand{\subm}{_{\mathrm{m}}}

\newcommand{\feas}{\mathcal{X}}
\newcommand{\pset}{\mathcal{P}_{\mathrm{S}}}
\newcommand{\pfront}{\mathcal{P}_{\mathrm{F}}}
\newcommand{\pcrit}{\mathcal{P}_{\mathrm{crit}}}

\newcommand{\gradf}{\bm{\nabla}f}
\newcommand{\jacobianf}{{\ve D{\ve f}}}

\newcommand{\gradm}{\bm{\nabla}m}

\newcommand{\gradmt}{\gradm\itat}

\newcommand{\hessianmt}[1]{\bm{H}{m_{#1}\ita{t}}}

\newcommand{\xitat}{\ve{x}\ita{t}}
\newcommand{\vemt}{\ve{m}\ita{t}}
\newcommand{\mitat}{m\itat}
\newcommand{\strictx}{\ve{{ \hat{x} }}_{\mathrm{PC}}\itat}
\newcommand{\stricts}{\ve{{ \hat{s} }}_{\mathrm{PC}}\itat}

\newcommand{\omga}[1]{ \omega\left( #1 \right) }
\newcommand{\omegamt}[1]{ \omega\subm\itat \left( #1 \right) }
\newcommand{\modomegamt}[1]{ \varpi\subm\itat \left( #1 \right) }
\newcommand{\modomegamtcrit}[2]{ \varpi_{\mathrm{m}_{ #1 }}\itat \left( #2 \right) }
\newcommand{\modomega}[1]{ \varpi\left( #1 \right) }

\newcommand{\modomegamit}[1]{ \varpi\subm\ita{#1} \left( \ve x\ita{#1} \right) }
\newcommand{\modomegait}[1]{ \varpi \left( \ve x\ita{#1} \right) }
\newcommand{\modshortmit}[1]{ \varpi\subm^{\left( #1 \right)} }
\newcommand{\modshortit}[1]{ \varpi^{\left( #1 \right)} }

\newcommand{\radius}{\Delta}
\newcommand{\radiust}{\radius\ita{t}}
\newcommand{\radiusmax}{{\radius}^{\mathrm{ub}}}

\newcommand{\radiusnext}{\radius\ita{t+1}}

\newcommand{\hessboundmt}{H_{\mathrm{m}}\ita t}
\newcommand{\hessboundm}{\mathrm{H}_{\mathrm{m}}}
\newcommand{\normconst}{\mathrm{c}}
\newcommand{\dmt}{\ve d\ita{t}_{\mathrm m}}
\newcommand{\norm}[1]{\left\| #1 \right\|}
\newcommand{\abs}[1]{\left| #1 \right|}
\newcommand{\normdmt}{\norm{\dmt}}
\newcommand{\xtrialt}{\ve x_+\ita t}

\newcommand{\stept}{\ve{s}\itat}
\newcommand{\xcauchyt}{\ve{x}_{\mathrm{PC}}\ita{t}}
\newcommand{\scauchyt}{\ve{s}_{\mathrm{PC}}\ita{t}}
\newcommand{\modscauchyt}{\ve{{ \tilde{s} }}_{\mathrm{PC}}\ita{t}}

\newcommand{\ball}[2]{B\left( #1 ; #2 \right)}
\newcommand{\ballt}{B\ita t}

\newcommand{\Phif}{\Phi}
\newcommand{\Phim}{\Phi_{\mathrm m}}
\newcommand{\Phimt}{\Phim\ita t}
\newcommand{\nusuccess}{\nu_{{\scriptscriptstyle++}}}
\newcommand{\nuaccept}{\nu_{\scriptscriptstyle+}}

\newcommand{\backtrackconst}{\mathrm{b}}
\newcommand{\armijofactor}{\mathrm{a}}
\newcommand{\sufficientconstant}{\kappa^{\scriptstyle\mathrm{sd}}}
\newcommand{\sufficientconstantm}{\sufficientconstant\subm}
\newcommand{\tildesufficientconstantm}{{\tilde{\kappa}}^{\mathrm{sd}}\subm}
\newcommand{\omegaconstant}{\kappa_{\omega}}

\newcommand{\rconstant}{\mathrm{r}}

\newcommand{\fullylinearconstant}{\epsilon}
\newcommand{\fullylinearconstantdf}{\dot{\fullylinearconstant}}

\newcommand{\successindices}{\mathcal{S}}

\newcommand{\lbomega}{\varpi_{\mathrm{m}}^{\mathrm{lb}}}% omega}}
\newcommand{\lbradius}{{\Delta}^{\mathrm{lb}}}

\newcommand{\gammasmallest}{\gamma_{\scriptscriptstyle\downdownarrows}}
\newcommand{\gammasmall}{\gamma_{\scriptscriptstyle\downarrow}}
\newcommand{\gammabig}{\gamma_{\scriptscriptstyle\uparrow}}
\newcommand{\epscrit}{\varepsilon_{\mathrm{crit}}}

\newcommand{\maximprovements}{\mathrm{M}}

\newcommand{\polyspacend}{\Pi_{n}^d}
\newcommand{\trainsites}{\Xi}

\newcommand{\maxiter}{\mathrm{N}_{\mathrm{it.}}}
\newcommand{\maxexpensive}{\mathrm{N}_{\mathrm{exp.}}}
\newcommand{\radiusmin}{\Delta_\mathrm{min}}
\newcommand{\radiuscrit}{\Delta_{\mathrm{crit}}}
\newcommand{\stepsmin}{\varepsilon_{\mathrm{rel}}}
\newcommand{\maxcritloops}{\mathrm{N}_{\mathrm{loops}}}
\newcommand{\ten}[1]{10^{#1}}

\begin{document}
%%%%%%%%%%%%%%%%%%%%%%%%%%%%%%%%%%%%%%%%%%

\section{Introduction}

Optimization problems arise in a multitude of applications in mathematics, 
computer science, engineering and the natural sciences.
In many real-life scenarios, there are multiple, equally important objectives
that need to be optimized. 
Such problems are then called \emph{Multiobjective Optimization Problems} (MOP).
In contrast to the single objective case, an MOP often does not have a single
solution but an entire set of optimal trade-offs between the different objectives,
which we call \emph{Pareto optimal}.
They constitute the \emph{Pareto Set} and their image is the \emph{Pareto Frontier}. 
The goal in the numerical treatment of an
MOP is to either approximate these sets or to find single points within these sets.
In applications, the problem can become more difficult when some of the objectives
require computationally expensive or time consuming evaluations.
For instance, the objectives could depend on a computer simulation or some
other \emph{black-box}.
It is then of primary interest to reduce the overall number of function evaluations.
Consequently, it becomes infeasible to approximate derivative 
information of the true objectives using, e.g., finite differences. 
In this work, optimization methods that do not use the objective gradients 
(which nonetheless are assumed to exist) are referred to as \emph{derivative-free}. 

There is a variety of methods to deal with multiobjective optimization problems, 
some of which are also derivative-free or try to constrain the number of expensive
function evaluations.
A broad overview of different problems and techniques concerning multiobjective optimization 
can be found, e.g., in \cite{ehrgott_multicriteria_2005,jahn,miettinen_nonlinear_2013,%
eichfelder_twenty_nodate}.
One popular approach for calculating Pareto optimal solutions is scalarization, 
i.e., the transformation of an MOP into a single objective problem, cf.\ \cite{eichfelder:adaptive} for an overview.
Alternatively, classical (single objective) descent algorithms can be adapted for
the multiobjective case \cite{fukuda_drummond,fliege_svaiter,grana_drummond_steepest_2005,%
nonlinear_cg,wolfe_linesearch,GPD19}. 
What is more, the structure of the Pareto Set can 
be exploited to find multiple solutions \cite{hillermeier_nonlinear_2001,GPD19b}.
There are also methods for non-smooth problems \cite{wilppu_new_2014,%
gebken_efficient_2021} and multiobjective direct-search variants 
\cite{direct_search,audet_multiobjective}.
Both scalarization and descent techniques may be included in 
Evolutionary Algorithms (EA)
\cite{kalyanmoydeb2001,coello_evo,abraham_evolutionary_2005,zitzler_diss},
the most prominent of which probably is NSGA-II \cite{deb_fast_2002}.
To address computationally expensive objectives or missing derivative information, 
there are algorithms that use surrogate models 
(see the surveys \cite{peitz_survey_2018,chugh_survey_2019,%
deb_surrogate_2020}) or borrow from ideas from scalar trust region methods, 
e.g., \cite{deb_trust_region_2019}.

In single objective optimization, trust region methods are well suited for derivative-free
optimization \cite{conn:derivative_free_2009,larson_wild_derivative_free}.
Our work is based on the recent development of multiobjective trust region methods:
\begin{itemize}
    \item In \cite{qu}, a trust region method using Newton steps for functions with positive 
    definite Hessians on an open domain is proposed.
    \item In \cite{villacorta} quadratic Taylor polynomials are used to compute the 
      steepest descent direction which is used in a backtracking manner to find solutions 
      for unconstrained problems.
    \item In \cite{ryu_derivative-free_2014} polynomial regression models are used 
      to solve an augmented MOP based on the scalarization in \cite{audet_multiobjective}.
      The algorithm is designed unconstrained bi-objective problems.
    \item In \cite{thomann}, quadratic Lagrange polynomials are used and 
    the Pascoletti-Serafini scalarization is employed for the descent step calculation.
\end{itemize}
Our contribution is the extension of the above-mentioned methods to general fully linear models (and in particular radial basis function surrogates as in \cite{wild:orbit}),
%into a flexible algorithmic framework for fully linear models, 
which is related to the scalar framework in \cite{conn:trm_framework}.
Most importantly, this reduces the complexity with respect to the parameter space dimension to linear, in contrast to the quadratically increasing number of function evaluations in other methods.
We further prove convergence to critical points when the problem is 
constrained to a convex and compact set by using an analogous 
argumentation as in \cite{conn_trust_region_2000}.
This requires new results concerning the continuity of the projected steepest descent direction.
We also show how to keep the convergence properties for constrained problems when the 
Pascoletti-Serafini scalarization is employed (like in \cite{thomann}).

The remainder of the paper is structured as follows: 
\cref{section:optimality} provides a brief introduction to multiobjective optimality and 
criticality concepts.
In \cref{section:algo_description} the fundamentals of our algorithm are explained.
In \cref{section:models} we introduce fully linear surrogate models and describe 
their construction.
We also formalize the main algorithm in this section.
\cref{section:steps} deals with the descent step calculation so that a sufficient decrease
is achieved in each iteration.
Convergence is proven in \cref{section:convergence} and a few numerical examples are 
shown in \cref{section:examples}. 
We conclude with a brief discussion in \cref{section:conclusion}.

\section{Optimality and Criticality in Multiobjective Optimization}
\label{section:optimality}

We consider the following (real-valued) multiobjective optimization problem:
\begin{equation}
	\min_{\ve x\in \feas} \ve f(\ve x) :=
	\min_{\ve x\in \feas}
	\begin{bmatrix}
		f_1(\ve x)\\
		\vdots\\
		f_k(\ve x)\\
	\end{bmatrix} \in \mathbb{R}^k,
	\tag{MOP}
	\label{eqn:mop}
\end{equation}
with a feasible set $\feas \subseteq \mathbb{R}^n$ and $k$ objective functions 
$f_\ell\colon \mathbb{R}^n \to \mathbb{R},\ \ell=1,\ldots,k$.
We further assume \eqref{eqn:mop} to be \emph{heterogeneous}.
That is,
there is a non-empty subset $\expensiveindices\subseteq \{1,\ldots,k\}$ of indices
so that the gradients of $f_\ell, \ell\in \expensiveindices,$ are unknown and cannot be approximated, 
e.g., via finite differences.
The (possibly empty) index set $\cheapinidces = \{1,\ldots,k\}\setminus \expensiveindices$ indicates 
functions whose gradients are available.

Solutions for \eqref{eqn:mop} consist of optimal trade-offs $\ve x^*\in \feas$ between the 
different objectives and are called non-dominated or Pareto optimal. 
That is, there is no $\ve x\in \feas$ with $\ve f(\ve x) \prec \ve f(\ve x^*)$ 
(i.e., $\ve f(\ve x) \le \ve f(\ve x^*)$ and $f_\ell(\ve x) < f_\ell (\ve x^*)$ 
for some index $\ell\in\{1,\ldots,k\}$). 
The subset $\pset \subseteq \feas$ of non-dominated points is then called the
\emph{Pareto Set} and its image $\pfront:= \ve f\left(\pset \right)\subseteq \mathbb{R}^k$ is 
called the \emph{Pareto Frontier}.
% We will call $\ve y \in \pfront$ \emph{efficient.} 
All concepts can be defined in a local fashion in an analogous way.

Similar to scalar optimization, local optima can be characterized using the gradients 
of the objective function. We therefore implicitly assume all 
objective functions $f_\ell, \ell=1,\ldots,k,$ to be continuously differentiable on $\feas$.
Moreover, the following assumption allows for an easier treatment of 
tangent cones in the constrained case: 
\begin{Assumption}
	Either $\feas = \mathbb{R}^n$ or the feasible set $\feas \subseteq \mathbb{R}^n$ is closed, bounded and convex.
	All functions are defined on $\feas$.
	\label{assumption:feasible_set_compact_convex}
\end{Assumption}

Because $\mathbb{R}^k$ is finite-dimensional \cref{assumption:feasible_set_compact_convex} 
is equivalent to requiring $\feas$ to be compact and convex, 
which is a standard assumption in the MO literature \cite{fliege_svaiter,fukuda_drummond}.\\
Now let $\gradf_\ell(\ve x)$ denote the gradient of $f_\ell$ and 
$\jacobianf(\ve x) \in \mathbb{R}^{k\times n}$ the Jacobian of $\ve f$ at $\ve x\in \feas$.

\begin{Definition}
	We call a vector $\ve d \in \feas - \ve x$ a multi-descent direction for $\ve f$ in 
	$\ve x$ if $\langle \gradf_\ell(\ve x), \ve d \rangle < 0 \quad$ for all $\ell\in\{1,\ldots,k\},$
	or equivalently if
	\begin{equation}
  	    \max_{\ell=1,\ldots,k} \langle \gradf_\ell(\ve x^*), \ve d \rangle < 0
		\label{eqn:multi_descent_max}
	\end{equation}
    where $\langle \bullet,\bullet\rangle$ is the standard inner product on $\mathbb{R}^n$ and 
    we consider $\feas - \ve x = \feas$ in the unconstrained case  $\feas = \mathbb{R}^n$.
\end{Definition}

A point $\ve x^*\in \feas$ is called \emph{critical} for \eqref{eqn:mop} iff there is no
$\ve d\in \feas - \ve x^*$ with \eqref{eqn:multi_descent_max}.
As all Pareto optimal points are also critical (cf.\ \cite{fukuda_drummond,luc} or \cite[Ch.\ 17]{jahn}), 
it is viable to search for optimal points by calculating points from the 
superset $\pcrit \supseteq \pset $ of critical points for \eqref{eqn:mop}.
One way to do so is by iteratively performing descent steps.
\citet{fliege_svaiter} propose several ways to compute suitable descent 
directions.
The minimizer $\ve d^*$ of the following problem is
known as the multiobjective steepest-descent direction. % of $\ve f$ in $\ve x\in \feas$.
\begin{equation}
	\begin{aligned}
        \min_{ \ve d \in \feas - \ve x }  \max_{\ell=1,\ldots,k}
            \langle \gradf_\ell(\ve x), \ve d \rangle
        \quad\text{s.t.}\quad \norm{\ve d} \le 1.
	\end{aligned}
	\tag{P1}
	\label{eqn:descent_direction_problem1}
\end{equation}
Problem \eqref{eqn:descent_direction_problem1} has an equivalent reformulation as
\begin{equation}
	\begin{aligned}
	\min_{ \ve d \in \feas - \ve x } \beta \qquad \text{s.t.}\quad
	\norm{\ve d} \le 1\quad \mbox{and}\quad
	\langle \gradf_\ell(\ve x), \ve d \rangle \le \beta ~\forall~\ell=1,\ldots,k,
	\end{aligned}
	\tag{P2}
	\label{eqn:descent_direction_problem2}
\end{equation}
which is a linear program, if $\feas$ is defined by linear constraints and the maximum-norm $\norm{ \bullet} =\norm{\bullet}_\infty$ is used \cite{fliege_svaiter}.
We thus stick with this choice because it facilitates implementation, but note that other choices are possible (see for example \cite{thomann}).

Motivated by the next theorem we can use the optimal value of either problem
as a measure of criticality, i.e., as a multiobjective pendant for the gradient norm.
As is standard in most multiobjective trust region works (cf. \cite{qu,villacorta,thomann}), 
we flip the sign so that the values are non-negative.
%The solutions of \eqref{eqn:descent_direction_problem1} or \eqref{eqn:descent_direction_problem2} have some desirable properties:

\begin{Theorem}
	\label{theorem:p1_properties}
    For $\ve x \in \feas$ let $\ve d^*(\ve x)$ be the minimizer of 
    \eqref{eqn:descent_direction_problem1} and $\omga{  \ve x  }$ be the negative optimal value, that is
	$$
	    \omga{\ve x} := - \max_{\ell=1,\ldots,k} \langle \gradf_\ell(\ve x), \ve d^*(\ve x) \rangle .
	$$
	Then the following statements hold:
	\begin{enumerate}
	    \item \label{p1_properties1} $\omga{  \ve x  } \ge 0$ for all $\ve x \in \feas$.
    	\item \label{p1_properties2} The function $\omega\colon \mathbb{R}^n \to \mathbb{R}$ is continuous.
	    \item \label{p1_properties3} The following statements are equivalent:
    	\begin{enumerate}
	        \item The point $\ve x\in \feas$ is \emph{not} critical.
	        \item $\omga{  \ve x  } > 0$.
	        \item $\ve d^*(\ve x) \ne \ve 0$.
	    \end{enumerate}
	\end{enumerate}
	Consequently, the point $\ve x$ is critical iff $\omga{  \ve x  } = 0$.
\end{Theorem}

\begin{proof}
For the unconstrained case all statements are proven in \cite[Lemma 3]{fliege_svaiter}. \\
The first and the third statement hold true for $\feas$ convex and compact by definition.
The continuity of $\omega$ can be shown similarly as in \cite{fukuda_drummond}, see \cref{appendix:omega_continuity}.
\end{proof}

With further conditions on $\ve f$ and $\feas$ the criticality measure $\omga{ \ve x }$ is even Lipschitz continuous 
and subsequently uniformly and Cauchy continuous:
%%%%%%%%%%%%%%%%

\begin{Theorem}
    \label{theorem:omega_uniformly_continuous}
    If $\gradf_\ell, \ell = 1, \ldots, k,$ are Lipschitz continuous and 
    \cref{assumption:feasible_set_compact_convex} holds, then the
    map $\omga{  \bullet  }$ as defined in \cref{theorem:p1_properties} is uniformly continuous.
\end{Theorem}

\begin{proof}
The proof for $\feas = \mathbb{R}^n$ is given by \citet{thomann_diss}.
A proof for the constrained case can be found in \cref{appendix:omega_continuity} as to 
not clutter this introductory section.
\end{proof}
%
%From \cref{theorem:omega_uniformly_continuous} we can deduce Cauchy continuity, 
%that is, from $\lim_{m\to \infty} \norm{ \ve x\ita{m} - \ve y\ita{m} } = 0 $
%it follows that
%$$
%\lim_{m\to \infty}
%	\left|
%		\omega\left( \ve x\ita{m}\right) - \omega\left( \ve y\ita{m}\right)
%	\right| = 0.
%$$
Together with \cref{theorem:p1_properties} this hints at $\omga{  \bullet  }$ being a criticality measure 
as defined for scalar trust region methods in \cite[Ch. 8]{conn_trust_region_2000}:

\begin{Definition}
	We call
	$
		\pi\colon \mathbb{N}_0\times \mathbb{R}^n\to \mathbb{R},
	$
	a criticality measure for \eqref{eqn:mop} if $\pi$ is Cauchy continuous with respect to its second argument and if
	$$
		\lim_{t\to \infty} \pi(t, \xitat) = 0
	$$
	implies that the sequence $\left\{\xitat\right\}$ asymptotically approaches a Pareto-critical 
	point.
	\label{mydef:criticality_measure}
\end{Definition}

\section{Trust Region Ideas}
\label{section:algo_description}

Multiobjective trust region algorithms closely follow the design of scalar approaches 
(see \cite{conn_trust_region_2000} for an extensive treatment).
Consequently, the requirements and convergence proofs in \cite{qu,villacorta,thomann}
for the unconstrained multiobjective case are fairly similar to those in \cite{conn_trust_region_2000}.
We will reexamine the core concepts to provide a clear understanding and point out the 
similarities to the scalar case.

%With the primary goal being the approximation of Pareto critical points of \eqref{eqn:mop} 
The main idea is to iteratively compute multi-descent steps $\stept$
in every iteration $t\in \mathbb{N}_0$.
We could, for example, use the steepest descent direction given by 
\eqref{eqn:descent_direction_problem1}.
%But because \eqref{eqn:descent_direction_problem1} requires knowledge of the objective gradients 
This would require knowledge of the objective gradients 
- which need not be available for objective functions with indices in $\expensiveindices$.
Hence, benevolent surrogate model functions
$$
	\ve m\ita{t} \colon \mathbb{R}^n \to \mathbb{R}^k, \; 
	\ve x \mapsto \ve m\ita{t}(\ve x) =
	\begin{bmatrix}
		m_1\ita{t}(\ve x),
		\ldots,
		m_k\ita{t}(\ve x)
	\end{bmatrix}^T,
$$
are employed.
Note, that for cheap objectives $f_\ell, \ell\in \cheapinidces,$ we could simply
use $m_\ell = f_\ell$ as long as these $f_\ell$ are twice continuously 
differentiable and have Hessians of bounded norm.

The surrogate models are constructed to be sufficiently accurate within a trust region
\begin{equation}
	\ballt := \ball{\xitat}{\radiust} =
	\left\{
		\ve x\in \feas: \norm{ \ve x - \ve x\ita t} \le \radiust
	\right\}, \quad \text{with $\norm{\bullet}=\norm{\bullet}_\infty$,}
	\label{eqn:def_ballt}
\end{equation}
around the current iterate $\ve x\ita{t}$.
The \emph{model} steepest descent direction $\dmt$ can then computed as the optimizer of 
the surrogate problem
\begin{equation}
	\begin{aligned}
	    \omegamt{ \ve x\ita{t} } := &-\min_{ \ve d \in \feas - \ve x } \beta \\
		\text{s.t. }&\norm{\ve d}\le 1, \text{ and }
		\langle \gradmt_\ell(\ve x), \ve d \rangle \le \beta \quad \forall\ell=1,\ldots,k.
	\end{aligned}
	\tag{Pm}
	\label{eqn:descent_direction_problemm}
\end{equation}

Now let $\sigma\itat>0$ be a step size.
The direction $\dmt$ need not be a descent direction for the true objectives 
$\ve f$ and the trial point $\xtrialt = \ve x\ita{t} + \sigma\ita{t} \dmt$ 
is only accepted if a measure $\rho\ita t$ of improvement and model quality
surpasses a positive threshold $\nuaccept$.
As in \cite{villacorta, thomann}, we scalarize the multiobjective problems by defining
\begin{equation*}
	\begin{aligned}
        \Phif(\ve x) 
        := 
			\max_{\ell=1,\ldots,k} f_\ell(\ve x),
		\qquad
        \Phimt(\ve x) 
    	:= 
            \max_{\ell=1,\ldots,k} m_\ell\ita{t}(\ve x).
	\end{aligned}
	%\label{eqn:Phi_def}
\end{equation*}
Whenever $\Phif(\ve x\ita{t}) - \Phif(\xtrialt) > 0$, there is a reduction in at least one objective function of 
$\ve f$ because of 
$$
	0 < \Phif(\ve x\ita{t}) - \Phif(\xtrialt) =
	  f_\ell(\ve x\ita{t}) - f_q(\xtrialt)
		\stackrel{\text{df.}}\le f_\ell(\ve x\ita{t}) - f_\ell(\xtrialt),
$$
where we denoted by $\ell$ the maximizing index in $\Phif(\ve x\itat)$ and 
by $q$ the maximizing index in $\Phif(\xtrialt)$.
\footnote{The abbreviation “$\text{df.}$” above the inequality symbol stands for “(by) definition” 
and is used throughout this document when appropriate.}
Of course, the same property holds for $\Phimt(\bullet)$ and $\vemt$.

Thus, the step size $\sigma\ita t > 0$ is chosen so that the step 
$\ve s\itat = \sigma\itat \dmt$ satisfies both $\ve x\ita{t} + \stept \in B\ita t$ 
and a “sufficient decrease condition” of the form 
$$ \Phimt( \xitat ) - \Phimt( \xitat + \stept )
    \ge
      \sufficientconstant
      \omga{\xitat}
      \min
      \left\{
        \mathrm{C} \cdot \omga{\xitat},
        1,
        \radiust
      \right\} \ge 0,
$$ 
with a constant $\mathrm C>0$, see \cref{section:steps}.
Such a condition is also required in the scalar case 
\cite{conn_trust_region_2000,conn:trm_framework}
and essential for the convergence proof in \cref{section:convergence}, 
where we show $\lim_{t\to \infty} \omga{\xitat} = 0$.

Due to the decrease condition the denominator in the ratio of actual versus predicted reduction,
\begin{equation}
	\rho\ita{t} =
	\begin{dcases}
		\frac{
			\Phif(\xitat) - \Phif(\xtrialt)
        }{
			\Phimt(\xitat) - \Phimt(\xtrialt)
        }	
        &\text{if $\xitat \ne \xtrialt$},\\
        0 
        &\text{if $\xitat = \xtrialt \Leftrightarrow \stept = \ve 0$,}
	\end{dcases}
	\label{eqn:rho_definition}
\end{equation}
is nonnegative.
A positive $\rho\itat$ implies a decrease in at least one objective $f_\ell$,
so we accept $\xtrialt$ as the next iterate if $\rho\itat > \nuaccept > 0$.
If $\rho\ita{t}$ is sufficiently large, say $\rho\ita{t}\ge \nusuccess >\nuaccept > 0$, 
the next trust region might have a larger radius $\radiusnext\ge \radiust$.
%because we deem the surrogate models to be very accurate.
If in contrast $\rho < \nusuccess$, the next trust region radius should be smaller and 
the surrogates improved.\\
This encompasses the case $\stept = \ve 0$, when the iterate $\xitat$ is critical for
\begin{equation}
	\min_{\ve x\in B\ita t} \vemt (\ve x) 
	% :=
	% \min_{\ve x\in B\ita t}
	% \begin{bmatrix}
	% 	m_1\ita t(\ve x)\\
	% 	\vdots\\
	% 	m_k\ita t(\ve x)\\
	% \end{bmatrix} 
	\, \in \mathbb{R}^k.
	\tag{MOP$\mathrm{m}$}
	\label{eqn:mopmt}
\end{equation}
Roughly speaking, we suppose that $\xitat$ is near a critical point 
for the original problem \eqref{eqn:mop} if $\vemt$ is sufficiently accurate.
%Hence, a further improvement on a smaller subsequent trust region is advisable.
If we truly are near a critical point, then the trust region radius will approach 0.
For further details concerning the acceptance ratio $\rho\ita{t}$, see \cite[Sec. 2.2]{thomann}.

\begin{Remark}
	We can modify $\rho\itat$ in \eqref{eqn:rho_definition} 
	to obtain a descent in all objectives, i.e., if $\xitat \ne \xtrialt$ we 
	test 
	$
		\rho\itat = \dfrac{f_\ell(\xitat) - f_\ell(\xtrialt)}{
			m\itat_\ell(\xitat) - m\itat_\ell(\xtrialt)
		} >\nuaccept
	$
	for all $\ell = 1,\ldots,k$. This is the \emph{strict} acceptance test.
\end{Remark}

\section{Surrogate Models and the Final Algorithm}
\label{section:models}
Until now, we have not discussed the actual choice of surrogate models used for $\vemt$.
As is shown in \cref{section:steps}, the models should be twice continuously differentiable 
with uniformly bounded hessians.
To prove convergence of our algorithm we have to impose 
further requirements on the (uniform) approximation qualities of the surrogates $\vemt$.
We can meet these requirements using so-called fully linear models.
Moreover, fully linear models intrinsically allow for
modifications of the basic trust region method that are aimed at reducing 
the total number of expensive objective evaluations.
Finally, we briefly recapitulate how radial basis functions and 
multivariate Lagrange polynomials can be made fully linear.

\subsection{Fully Linear Models}
\label{section:fully_linear}
\noindent Let us begin with the abstract definition of full linearity as 
given in \cite{conn:trm_framework,conn:derivative_free_2009}:
\begin{Definition}
  Let $\radiusmax >0$ be given and let $f\colon \mathbb{R}\to \mathbb{R}$ be a function 
  that is continuously differentiable in an open domain containing 
  $\feas$ and has a Lipschitz continuous gradient on $\feas$.
  A set of model functions 
  $\mathcal M = \left\{ m\colon \mathbb{R}^n \to \mathbb{R} \right\}\subseteq C^1(\mathbb{R}^n, \mathbb{R})$ 
  is called a \emph{fully linear} class of models if the following hold:
  \begin{enumerate}
    \item There are positive constants $\fullylinearconstant,\fullylinearconstantdf$ and $L_m$ such that for any given $\radius \in (0, \radiusmax)$ and for any $\ve x \in \feas$ there is a model function $m \in \mathcal M$ with Lipschitz continuous gradient and corresponding Lipschitz constant bounded by $L_m$ and such that
    \begin{itemize}
    \item the error between the gradient of the model and the gradient of the function satisfies
    $$
      \norm{
        \gradf(\ve{\upxi}) - \gradm(\ve{\upxi})
      }
      \le \fullylinearconstantdf \radius, \quad \forall \ve{\upxi}\in \ball{\ve x}{ \radius },
    $$
    \item the error between the model and the function satisfies
    $$
      \left|
        f(\ve{\upxi}) - m(\ve{\upxi})
      \right|
      \le \fullylinearconstant \radius^2, \quad \forall \ve{\upxi}\in \ball{\ve x}{\radius }.
    $$
    \end{itemize}

    \item For this class $\mathcal M$ there exists “model-improvement” algorithm that 
    -- in a finite, uniformly bounded (w.r.t.\ $\ve x$ and $\radius$) number of steps --
    can
    \begin{itemize}
      \item either establish that a given model $m\in \mathcal{M}$ is fully linear on $\ball{\ve x}{\radius}$
      \item or find a model $\tilde m$ that is fully linear on $\ball{\ve x}{\radius}$.
    \end{itemize}
  \end{enumerate}
  \label{mydef:fully_linear}
\end{Definition}

\begin{Remark}
  In the constrained case, we treat the constraints as \emph{hard}, 
  that is, we do not allow for evaluations of the true objectives outside $\feas$, 
  see the definition of $\ballt\subseteq\feas$ in \eqref{eqn:def_ballt}.
  We also ensure to only select training data in $\feas$ 
  during the construction of surrogate models.

  In the unconstrained case, the requirements in \cref{mydef:fully_linear} can be relaxed a bit, 
  at least when using the strict acceptance test with 
  $\ve f(\ve x\ita{T}) \le \ve f(\ve x\ita t)$ for all $T\ge t\ge 0$.
  We can then restrict ourselves to the set
  \begin{align*}
    \feas^\prime :=
      \bigcup_{\ve x\in L(\ve x\ita 0)} \ball{\ve x}{\radiusmax},
      \quad\text{where }
      L(\ve x\ita{0}):=\left\{ \ve x\in \mathbb{R}^n: \ve f(\ve x) \le \ve f(\ve x\ita{0})\right\}.
  \end{align*}
\end{Remark}

For the convergence analysis in \cref{section:convergence}, 
we cite \cite[Lemma 10.25]{conn:derivative_free_2009} concerning 
the approximation quality of fully linear models on enlarged trust regions:
\begin{Lemma}
  %[\cite[Lemma 20.25]{conn:derivative_free_2009}]
  For $\ve x \in \feas$ and $\Delta\le \radiusmax$ consider a function $f$ and a fully-linear model $m$ as in \cref{mydef:fully_linear} with constants $\fullylinearconstant,
\fullylinearconstantdf, L_m>0$.
  Let $L_f>0$ be a Lipschitz constant of $\gradf$.\\
  Assume w.l.o.g.\ that
  $$
    L_m + L_f \le \fullylinearconstant \quad
    \text{and}\quad
    \frac{\fullylinearconstantdf}{2} \le \fullylinearconstant.
  $$
  Then $m$ is fully linear on $\ball{\ve x}{\tilde \Delta}$ for any $\tilde \Delta\in [\Delta,\radiusmax]$ with respect to the same constants $\fullylinearconstant,\fullylinearconstantdf,L_m$.
  \label{theorem:fully_linear_in_larger_region}
\end{Lemma}

\subsubsection{Algorithm Modifications}
\label{section:algo_mods}

\newcommand{\criticalityRoutine}{\texttt{criticalityRoutine}}

With \cref{mydef:fully_linear} we have formalized our assumption that the surrogates become more accurate when we decrease the trust region radius.
This motivates the following modifications:
\begin{itemize}
  \item 
    “Relaxing” the (finite) surrogate construction process to try for a possible descent even if the surrogates are not fully linear.
  \item 
    A criticality test depending on $\modomegamt{\xitat}$.
    If this value is very small at the current iterate, 
    then $\xitat$ could lie near a Pareto-critical point.
    With the criticality test and \criticalityRoutine{} we ensure that the next model is 
    fully linear and the trust region 
    is not too large.
    This allows for a more accurate criticality measure and 
    descent step calculation.
  \item 
    A trust region update that also takes into 
    consideration $\modomegamt{\xitat}$.
    The radius should be enlarged if we have a large acceptance ratio 
    $\rho\itat$ and the $\radiust$ is small as 
    measured against $\beta \omegamt{\xitat}$ for a constant $\beta> 0$.
\end{itemize}
These changes are implemented in \cref{algo:algorithm1}.
For more detailed explanations we refer to \cite[Ch. 10]{conn:derivative_free_2009}.

\begin{algorithm}[!htb]
  \Configuration{Criticality parameters $\epscrit > 0$ and $\mu > \beta > 0$, 
    acceptance parameters $1>\nusuccess \ge \nuaccept \ge 0, \nusuccess\ne 0$, 
    update factors 
    $\gammabig \ge 1 > \gammasmall \ge \gammasmallest > 0$ 
    and $\radiusmax>0$;
  }
  \KwIn{The initial site $\ve x \ita{0}\in \mathbb{R}^n$;}
%Initialize: Set $t=0$ und construct models $\ve m\itat$\;
  \For{$t=0,1,\ldots$}{
    \eIf{$t>0$ and iteration $(t-1)$ was model-improving (cf. \cref{mydef:iter_categories})}
    {
      Perform at least one improvement step on $\ve m\ita{t-1}$
      and then let $\ve m\itat \leftarrow \ve m\ita{t-1}$\;
    }
    {
      Construct surrogate models $\vemt$ on $\ballt$\;
    }
    \tcc{Criticality Step:}
    \If{
      $\modomegamt{\xitat} < \epscrit$ {\upshape\textbf{and}
      $\bigl($
        $\vemt$ not fully linear 
        {\upshape\textbf{or}}
        $\radiust > \mu\modomegamt{\xitat}$
      $\bigr)$}
    }
    {
      Set $\radiust_* \leftarrow \radiust$\;
      Call $\text{\criticalityRoutine}()$ so that $\vemt$ is fully linear on $\ballt$ with $\radiust\in (0, \mu\modomegamt{\xitat}]$\;
      Then set $\radiust \leftarrow
      %\max\left\{\radiust, \beta \modomegamt{\xitat} \right\}
      \min\left\{ \max\left\{\radiust, \beta \modomegamt{\xitat} \right\}, \radiust_* \right\}
      $
      \;
      %\tcc*{still $\le \mu\modomegamt{\xitat}$}
    }
    Compute a suitable descent step $\stept$\;
    Set $\xtrialt \leftarrow \xitat + \stept$, evaluate $\ve f(\xtrialt)$ and compute $\rho\ita t$ with \eqref{eqn:rho_definition}\;
    \DontPrintSemicolon
    Perform the following updates:\;
    $$
    \begin{aligned}
      \ve x\ita{t+1} &\leftarrow
        \begin{cases}
          \xitat & \text{if $\rho\ita t < \nuaccept$},\\
          \xtrialt & \text{if $\nuaccept \le \rho\itat < \nusuccess$ \& $\vemt$ is fully linear,}\\
          \xtrialt & \text{if $\nusuccess \le \rho \itat$}.
        \end{cases}
      \\
      \radiusnext &\leftarrow \radius_+, \text{ where }\\
        \radius_+ & 
        \begin{cases}
          = \radiust
            &\text{if $\rho\ita t < \nusuccess$ \& $\vemt$ is \textbf{not} fully linear},\\
          \in [\gammasmallest \radiust, \gammasmall\radiust] 
            &\text{if $\rho\ita t < \nusuccess$ \& $\vemt$ is fully linear},\\
          \in \left[\radiust, \min\{\gammabig \radiust, \radiusmax\} \right] 
            &\text{if $\nusuccess \le \rho\ita t$ and $\radiust \ge \beta\omegamt{\xitat}$,}\\
          = \min\{\gammabig \radiust, \radiusmax\} 
            &\text{if $\nusuccess \le \rho\ita t$ and $\radiust < \beta\omegamt{\xitat}$.}\\
        \end{cases}
    \end{aligned}
    $$\;
  \PrintSemicolon
  }
  \caption{General TRM for \eqref{eqn:mop}}
  \label{algo:algorithm1}
\end{algorithm}

\begin{procedure}[htb]
  \Configuration{backtracking constant $\alpha\in(0,1)$, $\mu>\beta>0$ from \cref{algo:algorithm1};}
  %\If{$\vemt$ not fully linear on $\ballt$}
  Set $\Delta_0 \leftarrow \radiust$\;
  \For{ $j = 1, 2,\ldots $ }{
    Set radius: $\radiust\leftarrow \alpha^{j-1} \Delta_0$\;
    Make models $\vemt$ fully linear on $\ballt$
    \tcc*{can change $\modomegamt{\xitat}$}
    \If{ $\radiust \le \mu \modomegamt{\xitat}$ }{
      Break\;
    }
  }
\caption{criticalityRoutine($ $)}
\end{procedure}

From \cref{algo:algorithm1} we see that we can classify the iterations 
based on $\rho\itat$ in the following way:
\begin{Definition}
  For given constants $0\le \nuaccept\le \nusuccess < 1, \nusuccess\ne 0,$ 
  we call the iteration with index $t\in \mathbb{N}_0$ of 
  \cref{algo:algorithm1} \ldots
  \begin{itemize}
    \item \ldots \emph{successful} if $\rho\itat \ge \nusuccess$. The set of successful indices is
      $\successindices = \{ t\in \mathbb{N}_0: \rho\itat \ge \nusuccess\} \subseteq \mathbb{N}_0$.
    \item \ldots \emph{model-improving} if $\rho\itat < \nusuccess$ and the models 
      $\vemt = [m_1\itat, \ldots, m_k\itat]^T$ are not fully linear.
      In these iterations the trust region radius is not changed.
    \item \ldots \emph{acceptable} if $\nuaccept > \rho\itat \ge \nuaccept$ and 
      the models $\vemt$ are fully linear.
      If $\nusuccess = \nuaccept \in (0,1)$, then there are no acceptable indices.
    \item \ldots \emph{inacceptable} otherwise, 
      i.e., if $\rho\itat < \nuaccept$ and $\vemt$ are fully linear.
  \end{itemize}
  \label{mydef:iter_categories}
  % The set indices that are not inacceptable is
  % $\acceptableindices = \{ t\in \mathbb{N}_0: \rho\itat \ge \nuaccept \}$, 
  % a superset of $\successindices$.
\end{Definition}

% This classification is illustrated in \cref{fig:iter_categories}.
%%% THERE IS AN ERROR IN THE GRAPHIC :
%%% MODEL-IMPROVING ITERATIONS DO NOT HAVE TO ACCEPT XTRIALT
% \begin{figure}[H]
%   \begin{center}
%     \resizebox{!}{.5\textwidth}{
%       \input{figs/tex/step_classification.tex}
%     }
%   \end{center}
%   \caption{Classification of the iterations of based on $\rho\itat$ and the respective updates.}
%   \label{fig:iter_categories}
% \end{figure}

%%% (end) input of file  /home/manuelbb/Desktop/latex_cleanup/trm_combined/sections/5_models/fully_linear.tex  %%%%

\subsection{Fully Linear Lagrange Polynomials}
% The error bounds in \cref{mydef:fully_linear} are conceptually similar to the 
% error estimates obtained by constructing first degree Taylor polynomial models 
% for twice differentiable objectives using a central difference scheme.
% Now, it would of cause be desirable to induce the surrogates with
% higher order information of the true objectives.
% But because we assume expensive objectives, it is not feasible to construct 
% even second degree finite difference models which require 
% $\mathcal O(n^2)$ evaluations.\\
Quadratic Taylor polynomial models are used very frequently. 
As explained in \cite{conn:derivative_free_2009} we can alternatively use 
multivariate interpolating Lagrange polynomial models when derivative information 
is not available.
We will consider first and second degree Lagrange models.
Even though the latter require $\mathcal O(n^2)$ function evaluations 
% to build 
they are still cheaper than second degree finite difference 
models. 
% Whilst the latter also require $\mathcal O(n^2)$ evaluations,
% the factor is usually much lower than for second degree finite difference 
% models. 
For this reason, these models are also used in \cite{thomann,thomann_diss}.

% %\subsubsection{Poised Training Sets}
% The polynomial models belong to the space 
% $\polyspacend$ of real-valued $n$-variate polynomials of degree $d$.
% The dimension of $\polyspacend$ is 
% % $$
% %   p := \dim \polyspacend = \binom{n+d}{d} = \binom{n+d}{n},
% %   \quad\text{(see \cite[Th. 2.5]{wendland})}
% % $$
% $p = n+1$ for $d=1$ and $p=\dfrac{(n+1)(n+2)}{2}$ for $d=2$.
% In contrast to the univariate case, 
% % $n=1$ we cannot interpolate arbitrary data with 
% % models in $\polyspacend$.
% for $n\ge 2$ the \emph{Mairhuber-Curtis} theorem\cite{wendland} applies and the 
% interpolation sites must form a so-called $\polyspacend$-unisolvent or \emph{poised} 
% subset of order $p$ in $\feas$.
To construct an interpolating polynomial model we have to provide $p$ data 
sites, where $p$ is the dimension of the space
$\polyspacend$ of real-valued $n$-variate polynomials with degree $d$.
For $d=1$ we have $p = n+1$ and for $d=2$ it is $p=\dfrac{(n+1)(n+2)}{2}$.
If $n\ge 2$, the \emph{Mairhuber-Curtis} theorem\cite{wendland} applies and the 
data sites must form a so-called \emph{poised} set in $\feas$.
% subset of order $p$ in $\feas$.
%\begin{Definition}
  The set $\trainsites = \{ \ve{\upxi}_1, \ldots, \ve{\upxi}_p \} \subset\mathbb R^n$
  is poised 
% for interpolation in $\polyspacend$ if the matrix
  if for any basis $\{\psi_i\}$ of $\polyspacend$ the matrix
  % $
  %  \ve M_\psi := 
  %   \begin{bmatrix}
  %     \ve{\uppsi}(\ve{\upxi}_1) & \ldots & \ve{\uppsi}(\ve{\upxi}_p)
  %   \end{bmatrix}  \ \in \mathbb R^{p\times p}
  % $
  $
  \ve M_\psi := \left[ \psi_i(\ve{\upxi}_j) \right]_{1\le i,j\le p}
  $
  is non-singular.
% , where the columns evaluate the basis,
%  i.e., $\ve{\uppsi}(\ve{\upxi}) = [\psi_1(\ve{\upxi}), \ldots, \psi_p(\ve{\upxi})]^T$.
%\end{Definition}
Then there is a unique polynomial 
$m(\ve{x}) = \sum_{i=1}^p \lambda_i \psi_i(\ve{x})$ with 
$m(\ve{\upxi}_j) = F(\ve{\upxi}_j)$ for all $j=1,\ldots,p$ and any function $F\colon \mathbb R^n \to \mathbb R$.
Given a poised set $\trainsites$ the associated Lagrange basis $\{l_i\}$ 
of $\polyspacend$ is defined by $l_i(\ve{\upxi_j}) = \delta_{i,j}$.
The model coefficients then simply are the data values, 
i.e., $\lambda_i = F(\ve{\upxi}_i)$.

Same as in \cite{thomann_diss}, we implement Algorithm 6.2 from 
\cite{conn:derivative_free_2009} to ensure poisedness.
It selects training sites $\trainsites$ from the current (slightly enlarged)
trust region of radius $\theta_1 \radiust$ and calculates the associated 
lagrange basis.
We can then separately evaluate the true 
objectives $f_\ell$ on $\trainsites$ to easily build 
the surrogates $m_\ell\itat$, $\ell \in \{1,\ldots, k\}$.
Our implementation always includes $\ve{\upxi}_1 = \xitat$ and 
tries to select points from a database of prior evaluations first.

% In \cite{thomann,thomann_diss}, it is stated as an assumption 
% that these models satisfy
% the error bounds in \cref{mydef:fully_linear}.
%We get rid of this requirement 
We employ an additional algorithm 
(Algorithm 6.3 in \cite{conn:derivative_free_2009}) to ensure that the set 
$\trainsites$ is even \emph{$\Lambda$-poised}, 
see \cite[Definition 3.6]{conn:derivative_free_2009}.
The procedure is still finite and ensures 
%that $\trainsites$ spans well $\ve{\uppsi}(\ballt)$ and 
the models are actually \emph{fully linear}.
The quality of the surrogate models can be improved by choosing a 
small algorithm parameter $\Lambda > 1$.
Our implementation tries again to recycle points from a database. 
Different to before, interpolation at $\xitat$ can no longer be guaranteed.
This second step can also be omitted first and then used as a model-improvement  
step in a subsequent iteration.

% To wrap it up, the model construction step in \cref{algo:algorithm1} consists
% of the following steps for fully linear Lagrange models:
% \begin{enumerate}
%   \item Determine a poised set $\trainsites = \{ \ve{\upxi}_1, \ldots, 
%     \ve{\upxi}_p \}$ within the current (possibly enlarged) trust region $\ballt$.
%   \item Make the set $\trainsites$ be $\Lambda$-poised (for a fixed $\Lambda>1$).
%   \item Calculate the Lagrange basis $\{l_1, \ldots, l_p\}$ associated with
%     $\trainsites$. 
%   \item For each $\ell = 1, \ldots , k,$ set 
%   $$
%     m_\ell\itat(\ve x) = \sum_{i=1}^p f(\ve{\upxi}_i) l_i(\ve x).
%   $$
% \end{enumerate}
% One could also omit step 2 first, try a descent step and perform step 2 as a 
% model-improvement step in the next iteration.

\subsection{Fully Linear Radial Basis Function Models}
\label{section:rbf}
The main drawback of quadratic Lagrange models is that we still need $\mathcal O(n^2)$
function evaluations in each iteration of \cref{algo:algorithm1}.
A possible fix is to use under-determined regression polynomials instead 
\cite{conn:derivative_free_2009,wild:dissertation,ryu_derivative-free_2014}.
Motivated by the findings in \cite{wild:orbit} we chose so-called Radial Basis Function (RBF)
models as an alternative.
RBF are well-known for their approximation capabilities on irregular data \cite{wendland}.
In our implementation they have the form 
\begin{equation}
  m(\ve x) = \sum_{i=1}^N c_i \varphi\left( \norm{ \ve x - \ve{\upxi}_i }_2 \right) +  \pi (\ve x), 
  \quad \pi = \sum_{i} \lambda_i \psi_i \in \polyspacend, d\in \{0,1\}, N \in \mathbb N, 
  \label{eqn:rbf_form}
\end{equation}
where 
% $\pi = \sum_{i} \lambda_i \psi_i \in \polyspacend$ is a polynomial with 
% $d = \deg \pi \in \{0,1\}$ and 
$\varphi$ is a function from $\mathbb R_{\ge 0}$ to $\mathbb R$.
For a fixed $\varphi$ the mapping $\varphi\left( \norm{\bullet} \right)$ from 
$\mathbb R^n \to \mathbb R$ is radially symmetric with respect to its argument
and the mapping 
$(\ve x, \ve{\upxi} )\mapsto \varphi\left( \norm{ \ve x - \ve{\upxi} }_2 \right)$ 
is called a \emph{kernel}.

\citet{wild:orbit} describe a construction of RBF surrogate models as in 
\eqref{eqn:rbf_form} (see also \cite{wild_global_2011} and the 
dissertation \cite{wild:dissertation} for more details). 
If we restrict ourselves to functions 
$\varphi\left( \norm{\bullet} \right)$
that are conditionally positive definite 
(c.p.d. -- see \cite{wendland,wild:orbit} for the 
definition) of order at most two, then the surrogates can be made 
certifiably fully linear with $N = n + 1$.
As before, the algorithms tries to select an initial training set 
$\trainsites = \{ \ve{\upxi}_1, \ldots, \ve{\upxi}_N\}\subset B(\xitat; \theta_1 \radiust)$
with $N = n+1$ and a scaling factor $\theta_1\ge 1$. 
The set must be poised for interpolation with affine linear polynomials.
Due to $\varphi\left( \norm{\bullet} \right)$ being c.p.d. of 
order $ D \le 2$, the interpolation system 
$$
\begin{aligned}
&\begin{bmatrix}
  \ve{\Upphi} & \ve M_\psi^T \\
  \ve M_\psi & \ve 0
\end{bmatrix}
\cdot 
\begin{bmatrix}
  \ve c \\ \ve \uplambda 
\end{bmatrix} = \begin{bmatrix}
  \ve F( \trainsites ) \\ \ve 0 
\end{bmatrix}  
  \in \mathbb R^{N + p \times 1},
  \\
  &
  \ve c = [ c_i ]_{1\le i \le N}, \
  \ve{\uplambda} = [ \lambda_i ]_{1\le i \le p}, \
  \ve{\Upphi} = [\varphi\left( \norm{ \ve{\upxi}_i - \ve{\upxi}_j } \right) ]_{1 \le i,j \le N},
\end{aligned}
$$
is uniquely solvable for any $F\colon \mathbb R^n \to \mathbb R$ if we 
choose $\polyspacend$ such that $d \ge \max\{ 0, D - 1\}$.
We can even include more points, $N\ge n+1$, from within a region of 
maximum radius $\theta_2\radiusmax$, $\theta_2\ge \theta_1\ge 1$,
to capture nonlinear behavior of $F$.
% Without going too much into detail, the model construction step in \cref{algo:algorithm1}
% for RBF surrogates is as follows, with enlargement factors $\theta_2 \ge \theta_1 \ge 1$:
% \begin{enumerate}
%   \item Try to find a poised set $\trainsites$ with $N = n+1$ sites from a database 
%     within the enlarged trust region of radius $\theta_1 \radiust$. 
%   \item Let $|\trainsites| = N^*$.
%     If $N^* = n+1$, then the models will be fully linear. Go to step 4.\\
%     Else try to complete $\trainsites$ to be poised with points within a 
%     box of radius $\theta_2 \radiusmax$. 
%     If now $|\trainsites| = n+1$ go to step 4.
%   \item If it still holds that $|\trainsites| = N^*$, then the model will be fully linear.\\
%     Evaluate all objectives at $n+1 - |\trainsites|$ suitable \textbf{new} sites
%     within $\theta_1\radiust$.
%   \item Try to include up to more points from a database, so that $|\trainsites| \le N_{\text{max}}$.
%   \item For every $\ell = 1, \ldots, k,$ solve the interpolation system 
%   $\{ m_\ell\itat(\ve{\upxi}) = f_\ell(\ve{\upxi}) \}_{\ve{\upxi}\in \trainsites}$.
% \end{enumerate}
More detailed explanations can be found in \cite{wild:orbit}. 
Modifications for box constraints are shown in \cite{wild:dissertation} and \cite{conorbit}.\\

\begin{specialtable}[H] 
  \caption{
    Some radial functions $\varphi\colon \mathbb R_{\ge 0}\to \mathbb R$ 
    that are c.p.d. of order $D\le 2$, cf. \cite{wild:orbit}.%
  }
  \label{table:rbf}
  %%% \tablesize{} %% You can specify the fontsize here, e.g., \tablesize{\footnotesize}. If commented out \small will be used.
  \begin{tabular}{cccc}
  \toprule
  \textbf{Name} & $\varphi(r)$  & c.p.d. order $D$	& $\deg \pi$ \\
  \midrule
  Cubic	        & $r^3$ & $2$ & $1$ \\
  Multiquadric	& $-\sqrt{1 + (\alpha r)^2}, \alpha > 0 $  & $1$ & $\{0,1\}$ \\
  Gaussian      & $\exp( - (\alpha r)^2 ), \alpha > 0 $    & $0$ & $\{0,1\}$ \\
  \bottomrule
  \end{tabular}
\end{specialtable}

\begin{figure}[H]
  \centering
  \includegraphics[width=\linewidth]{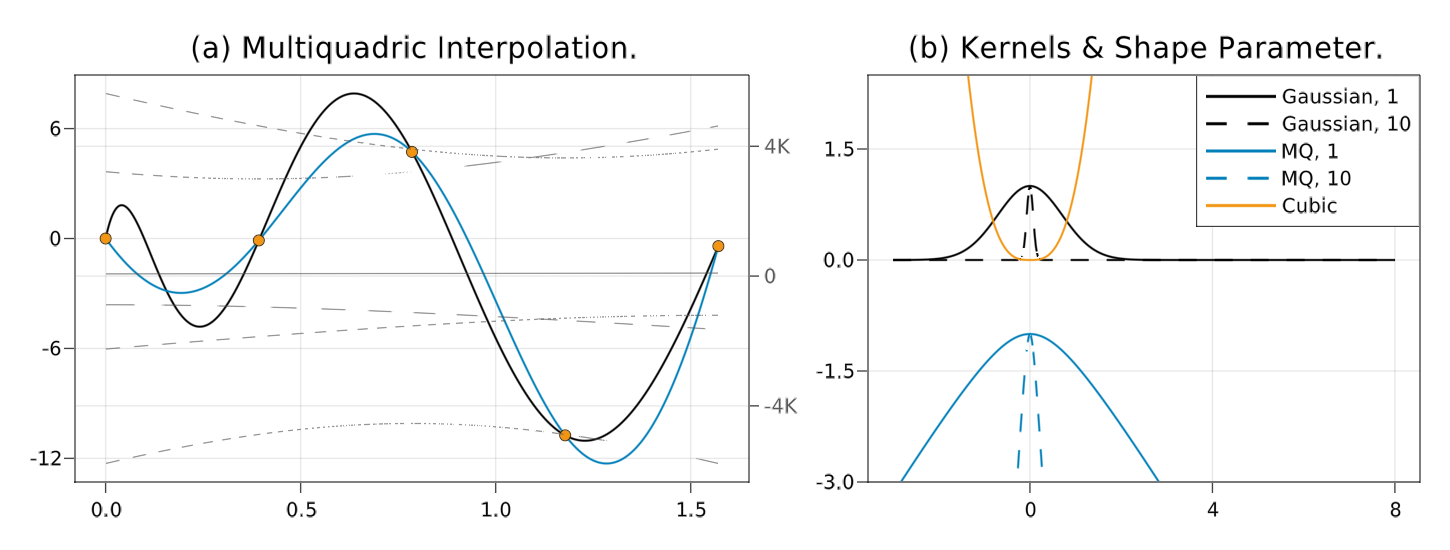}\\
  \caption{%
    \textbf{(a)} Interpolation of a nonlinear function (black) by a Multiquadric surrogate (black) based on 5 discrete training points (orange).
    Dashed lines show the kernels and the polynomial tail.
    \textbf{(b)} Different kernels in 1D with varying shape parameter ($1$ or $10$), see also \cref{table:rbf}.%
  }
  \label{fig:shapeparam}
\end{figure}

\cref{table:rbf} shows the RBF we are using and the possible polynomial degrees 
for $\pi$. 
Both the Gaussian and the Multiquadric allow for fine-tuning with a 
shape parameter $\alpha>0$.
This can potentially improve the conditioning of the interpolation system.
\cref{fig:shapeparam}~\textbf{(b)} illustrates the effect of the shape parameter.
As can be seen, the radial functions become narrower for larger shape parameters.
Hence, we do not only use a constant shape parameter $\alpha=1$ like \citet{wild:orbit} do,
but we also use an $\alpha$ that is (within lower and upper bounds) inversely 
proportional to $\radiust$.\\
\cref{fig:shapeparam}~\textbf{(a)} shows interpolation of a nonlinear function by 
a surrogate based on the Multiquadric with a linear tail.

\section{Descent Steps}
\label{section:steps}
In this section we introduce some possible steps $\stept$ to use in \cref{algo:algorithm1}.
We begin by defining the best step along the steepest descent direction 
as given by \eqref{eqn:descent_direction_problemm}.
Subsequently, backtracking variants are defined that use a multiobjective variant of Armijo's 
rule.

\subsection{Pareto-Cauchy Step}
Both the \emph{Pareto-Cauchy point} as well as a backtracking variant, 
the \emph{modified Pareto-Cauchy point}, are points along the 
descent direction $\dmt$ within $B\ita t$ so that a sufficient decrease 
measured by $\Phimt(\bullet)$ and $\omegamt{ \bullet }$ is achieved.
Under mild assumptions we can then derive a decrease in terms of $\omga{ \bullet }$.

\begin{Definition}
	For $t\in \mathbb{N}_0$ let $\dmt$ be a minimizer for \eqref{eqn:descent_direction_problemm}.
  The best attainable trial point $\xcauchyt$ along $\dmt$ is called the 
  \emph{Pareto-Cauchy point} and given by
	\begin{equation}
		\begin{aligned}
		\xcauchyt &:= \xitat + \sigma\itat \cdot \dmt, \\
		&\sigma\itat = \argmin_{0\le \sigma} \Phimt
		\left(
			\xitat + \sigma \cdot \dmt
		\right)	&\text{s.t. $\xcauchyt\in B\ita{t}$.}
		\end{aligned}
		\label{eqn:pareto_cauchy_point}
	\end{equation}
	Let $\sigma\itat$ be the minimizer in \eqref{eqn:pareto_cauchy_point}.
	We call
	$\scauchyt := \sigma\itat\dmt$
	the \emph{Pareto-Cauchy step}.
	\label{mydef:true_pareto_cauchy}
\end{Definition}

If we make the following standard assumption, 
then the Pareto-Cauchy point allows for a lower bound on the improvement in terms of $\Phimt$.
%\newcounter{tmpcounter}
%\setcounter{tmpcounter}{\value{Assumption}}
%\setcounter{Assumption}{\value{diffassumption}}
\begin{Assumption}
    For all $t\in \mathbb{N}_0$ the surrogates $\ve m\itat(\ve x)= [m_1\itat(\ve x), \ldots, m_k\itat(\ve x)]^T$ 
    are twice continuously differentiable on an open set containing $\feas$.
    Denote by 
    %$\hessianf{\ell}(\ve x)$ the Hessian of
    % $f_\ell$ at $\ve x$ and by
    $\hessianmt{\ell}(\ve x)$ the Hessian of $m_\ell\itat$ for $\ell = 1,\ldots, k$.
	\label{assumption:f_m_twice_continuously_differentiable}
\end{Assumption}
%\setcounter{Assumption}{\value{tmpcounter}}
%\refstepcounter{Assumption}

\begin{Theorem}
  If 
  \cref{assumption:feasible_set_compact_convex,assumption:f_m_twice_continuously_differentiable} 
  are satisfied, then for any iterate $\xitat$ the Pareto-Cauchy point $\xcauchyt$ satisfies
	\begin{equation}
		\Phimt( \xitat ) - \Phimt( \xcauchyt ) \ge
		\frac{1}{2} \omegamt{  \xitat  } \cdot
		\min
		\left\{
			\frac{ \omegamt{ \xitat } }{ \normconst \hessboundmt },
			\radiust,
			1
		\right\},
		\label{eqn:sufficient_decrease_true_cauchy_point}
	\end{equation}
	where
	\begin{equation}
		\hessboundmt = \max_{\ell = 1, \ldots, k} \max_{\ve x\in \ballt} \norm{ \hessianmt{\ell}(\ve x) }_F
		\label{eqn:hessian_upper_bounds_ita}
	\end{equation}
    and the constant $\normconst>0$ relates the trust region norm $\norm{\bullet}$
    to the Euclidean norm $\norm{\bullet}_2$ via
	\begin{equation}
		\norm{ \ve x }_2 \le \sqrt{\normconst} \norm{ \ve x } \qquad \forall\ve x\in \mathbb{R}^n.
	\label{eqn:normconst}
	\end{equation}
\label{theorem:sufficient_decrease_true_cauchy_point}
\end{Theorem}

If $ \norm{ \bullet } = \norm{ \bullet}_\infty$ is used, 
then $\normconst$ can be chosen as $\normconst = k$.
The proof for \cref{theorem:sufficient_decrease_true_cauchy_point} is 
%in essence that given in \cite{villacorta} and 
provided after the next auxiliary lemma.

\begin{Lemma}
    \label{lemma:sufficient_decrease_true_cauchy_point}
    Under \cref{assumption:feasible_set_compact_convex,assumption:f_m_twice_continuously_differentiable}, 
    let $\ve d$ be a non-increasing direction at $\xitat \in \mathbb{R}^n$ for $\vemt$, i.e.,
    $$
    \left\langle \gradmt_\ell(\xitat), \ve d \right\rangle \le 0\qquad \forall\ell=1,\ldots, k.
    $$
    Let $q\in \{1,\ldots,k\}$ be any objective index and 
    $\bar \sigma \ge \min\left\{ \radiust, \norm{\ve d}\right\}$.
	Then it holds that
	\begin{equation*}
		m_q\itat(\xitat)
		- \min_{0<\sigma<\bar{\sigma}} m_q\itat
		\left( \xitat + \sigma \frac{{\ve d}}{\norm{{\ve d}}} \right)
		\ge
		\frac{w}{2}
	  \min \left\{
	   \frac{w}{\norm{  {\ve d}  }^2\normconst\hessboundmt},
	   \frac{\radiust}{\norm{  {\ve d}  }}
	 	,1
	  \right\},
  \end{equation*}
	where we have used the shorthand notation
	$$
		w = -\max_{\ell=1,\ldots,k}\left\langle
		\gradmt_\ell(\xitat), \ve d
		\right\rangle \ge 0.
	$$
	\label{theorem:minimize_along_direction}
\end{Lemma}

\cref{lemma:sufficient_decrease_true_cauchy_point} states that a minimizer 
along any non-increasing direction $\ve d$ achieves a minimum reduction
w.r.t.\ $\Phimt$.
Similar results can be found in in \cite{villacorta} or \cite{thomann}.
But since we do not use polynomial surrogates $\vemt$,
we have to employ the multivariate version of Taylor's theorem
to make the proof work.
%using the Lagrange form of the remainder.
We can do this because according to \cref{assumption:f_m_twice_continuously_differentiable}, the 
functions $m_q\ita t, q\in \{1,\ldots,k\}$ are twice continuously differentiable
in an open domain containing $\feas$.
Moreover, \cref{assumption:feasible_set_compact_convex} ensures that the function 
is defined on the line from $\ve{\upchi}$ to $\ve x$.
As shown in \cite[Ch. 3]{fleming} a first degree expansion at $\ve x \in B(\ve{\upchi}, \Delta)$ 
around $\ve{\upchi}\in \feas$ then leads to
\begin{equation}
\begin{aligned}
    &m_q\ita t(\ve x) = m_q( \ve{\upchi})
        + \gradmt_q(\ve{\upchi})^T \ve h
    + \frac{1}{2}\ve h^T \hessianmt{q}(\ve{\upxi}_q) \ve h, 
        & \text{with $\ve h = \left(\ve x -\ve{\upchi}\right)$},
    \\
    &\text{for some $\ve{\upxi}_q\in 
        \left\{ 
            \ve x + \theta\left(\ve{\upchi}-\ve x\right) : \theta\in [0,1] 
        \right\}$, for all $q=1,\ldots,k$.}
\end{aligned}
\label{eqn:multi_taylor}
\end{equation}

\begin{proof}[Proof of \cref{lemma:sufficient_decrease_true_cauchy_point}]
    \allowdisplaybreaks
    Let the requirements of \cref{lemma:sufficient_decrease_true_cauchy_point} hold and let $\ve d$ be a
    non-increasing direction for $\vemt$. Then:
\begin{align*}
    %&\vspace{-3em} \Phimt( \xitat ) - \Phimt( \xcauchyt )\\
    %& \stackrel{\mathclap{\eqref{eqn:Phi_def}}}\ge
    %m_q\itat(\xitat) - m_q\itat(\xcauchyt) \\
    %& \stackrel{ \mathclap{\eqref{eqn:pareto_cauchy_point}} } =
    &    m_q\itat(\xitat)
        - \min_{0<\sigma<\bar{\sigma}} m_q\itat
        \left( \xitat + \sigma \frac{{\ve d}}{\norm{{\ve d}}} \right)
    %\nonumber\\
    %&
    =
        \max_{ 0 \le \sigma \le \bar \sigma}
        \left\{
            m_q\itat(\xitat)
            - m_q\itat
                 \left( \xitat + \sigma \frac{{\ve d}}{\norm{{\ve d}}} \right)
        \right\}
    \nonumber\\
    %&\hspace{-4cm}
    &\stackrel{\mathclap{\eqref{eqn:multi_taylor}}}=
            \max_{ 0 \le \sigma \le \bar \sigma}
            \left\{
                m_q\itat(\xitat)
                - \left(
                    m_q\itat(\xitat)
                    + \frac{\sigma}{\norm{  {\ve d}  }} \langle \gradmt_q(\xitat), {\ve d}\rangle
                    + \frac{\sigma^2}{2\norm{  {\ve d}  }^2} \langle {\ve d}, \hessianmt{q} (\ve{\upxi}_q) {\ve d} \rangle
                \right)
            \right\}
    \nonumber\\
    % &=
    %     \max_{ 0 \le \sigma \le \bar \sigma}
    %     \left\{
    %         - \frac{\sigma}{\norm{  {\ve d}  }} \langle \gradmt_q(\xitat), {\ve d}\rangle
    %         - \frac{\sigma^2}{2\norm{  {\ve d}  }^2} \langle {\ve d}, \hessianmt{q} (\ve{\upxi}) {\ve d} \rangle
    %         \right\}
    % \nonumber\\
    &\ge
        \max_{ 0 \le \sigma \le \bar \sigma}
        \left\{
            - \frac{\sigma}{\norm{  {\ve d}  }}
            \max_{j=1,\ldots,k} \langle \gradmt_j(\xitat), {\ve d}\rangle
            - \frac{\sigma^2}{2\norm{  {\ve d}  }^2} \langle {\ve d}, \hessianmt{q} (\ve{\upxi}_q) {\ve d} \rangle
        \right\}.
    \shortintertext{
        We use the shorthand $w = - \max_j \langle \gradmt_j(\xitat), \ve d \rangle$ 
        and the Cauchy-Schwartz inequality to get}
    \dotso 
        &\ge
    \max_{ 0 \le \sigma \le \bar \sigma}
    \left\{
        \frac{\sigma}{\norm{  {\ve d}  }}
        w
        - \frac{\sigma^2}{2\norm{ \ve d }^2}
         \norm{ \ve d}_2^2 \norm{ \hessianmt{q} (\ve{\upxi})}_F
    \right\}
    %\nonumber\\
    %&
    \stackrel{
        %\mathclap{
            \eqref{eqn:normconst}, \eqref{eqn:hessian_upper_bounds_ita}
        %}
    }\ge
    \max_{ 0 \le \sigma \le \bar \sigma}
    \left\{
         \frac{\sigma}{\norm{  {\ve d}  }}
         w
         - \frac{\sigma^2}{2}\normconst \hessboundmt
    \right\}.
    \end{align*}
    The RHS is concave and we can thus easily determine the global maximizer $\sigma^*$.
    Similar to \cite[Lemma 4.1]{villacorta} we find
    \begin{align*}
     %\Phimt( \xitat ) - \Phimt( \xcauchyt )
     m_q\itat(\xitat)
     - \min_{0<\sigma<\bar{\sigma}} m_q\itat
      \left( \xitat + \sigma \frac{{\ve d}}{\norm{{\ve d}}} \right)
     &\ge
    %  \min \left\{
    %      \frac{
    %              w^2
    %          }{
    %              2\norm{  {\ve d}  }^2 \normconst \hessboundmt
    %          },
    %         \frac{\bar \sigma}{2\norm{  {\ve d}  }} w
    % \right\}\\
    % &=
    % \frac{w}{2}
    % \min \left\{
    %     \frac{w}{\norm{  {\ve d}  }^2\normconst\hessboundmt},
    %     \frac{\bar \sigma}{\norm{ {\ve d}}}
    % \right\}
    %  \\
    %  &\ge
    %  \frac{w}{2}
    %  \min \left\{
    %      \frac{w}{\norm{  {\ve d}  }^2\normconst\hessboundmt},
    %      \frac{\radiust}{\norm{\ve d}},
    %     \frac{\norm{\ve d}}{ \norm{  {\ve d}  }}
    %  \right\}
    %  \\
    % &=
     \frac{w}{2}
     \min \left\{
      \frac{w}{\norm{  {\ve d}  }^2\normconst\hessboundmt},
      \frac{\radiust}{\norm{  {\ve d}  }}
        ,1
     \right\},
    \end{align*}
    where we have additionally used $\bar\sigma \ge \min\{ \radiust, 1 \}$.
\end{proof}

\begin{proof}[Proof of \cref{theorem:sufficient_decrease_true_cauchy_point}]
If $\xitat$ is Pareto-critical for \eqref{eqn:mopmt}, then $\dmt = \ve 0$ and 
$\omegamt{\xitat} = 0$ and the inequality holds trivially.

Else, let the indices $\ell,q \in \{1, \ldots, k\}$ be such that
$$
  \Phimt( \xitat ) - \Phimt( \xcauchyt ) =
  m_\ell\itat (\xitat ) - m_q\itat (\xcauchyt) \ge m_q (\xitat) - m_q(\xcauchyt)
$$
and define
\begin{equation}
  \bar \sigma :=
    \begin{cases}
      \min \left\{ \radiust, \norm{ \dmt }\right\} 	
        &\text{if $\norm{\dmt} < 1$ or $\radiust \le 1$},\\
      \radiust
        &\text{else.}
    \end{cases}
  \label{eqn:bar_sigma_def}
\end{equation}
Then clearly $\bar \sigma \ge \min\left\{ \radiust, \norm{\dmt} \right\}$ and 
for the Pareto-Cauchy point we have
$$
    m_q\itat\left(\xcauchyt\right) 
    = 
        \min_{0\le \sigma\le \bar \sigma } m_q 
        \left( 
            \xitat + \frac{\sigma}{\norm{\dmt}} \dmt
        \right).
$$
From \cref{theorem:minimize_along_direction} and $\norm{\dmt}$ 
% it immediately follows that
% \begin{align*}
%   \Phimt( \xitat ) - \Phimt( \xcauchyt )
%     &\ge
%       \frac{ \omegamt{ \xitat } }{2}
%       \min \left\{
%         \frac{\omegamt{ \xitat }}{\norm{  {\ve d}  }^2\normconst\hessboundmt},
%         \frac{\radiust}{\norm{  {\ve d}  }}
%         ,1
%       \right\}
%     \shortintertext{and because of $\norm{\dmt}\le 1$, it follows that}
%     \Phimt( \xitat ) - \Phimt( \xcauchyt )
%     &\ge
%       \frac{ \omegamt{ \xitat } }{2}
%       \min
%       \left\{
%         \frac{\omegamt{ \xitat }}{\normconst\hessboundmt},
%         \radiust
%         ,1
%       \right\},
% \end{align*}
% which is what we needed to show.
the bound $\eqref{eqn:sufficient_decrease_true_cauchy_point}$ immediately follows.
\end{proof}

\begin{Remark}
	Some authors define the Pareto-Cauchy point as the actual minimizer $\ve x_{\min}\ita{t}$ of $\Phimt$ within the current trust region (instead of the minimizer along the steepest descent direction).
	For this true minimizer the same bound \eqref{eqn:sufficient_decrease_true_cauchy_point} holds.
    This is 
    %\spremark{
        due to
    %}
	\begin{align*}
	\Phimt(\xitat) - \Phimt( \ve x_{\min}\ita{t} )
	=
		m_\ell(\xitat) - \min_{\ve x\in \ballt} m_q(\ve x)
	\ge
		m_q(\xitat) - m_q(\xcauchyt).
	\end{align*}
	\label{remark:minimizer}
\end{Remark}

%%% (end) input of file  /home/manuelbb/Desktop/latex_cleanup/trm_combined/sections/4_steps/pareto_cauchy.tex  %%%%

%%% (start) input of file  /home/manuelbb/Desktop/latex_cleanup/trm_combined/sections/4_steps/modified_pareto_cauchy.tex  %%%%

\subsection{Modified Pareto-Cauchy Point via Backtracking}
A common approach in trust region methods is to find an approximate 
solution to \eqref{eqn:pareto_cauchy_point} within the current trust region.
Usually a backtracking approach similar to Armijo's inexact line-search 
is used for the Pareto-Cauchy subproblem.
Doing so, we can still guarantee a sufficient decrease.
%(as measured by $\Phimt$).
% \spout{that is a fraction of the decrease attained by the true Pareto-Cauchy step $\scauchyt$.}

Before we actually define the backtracking step along $\dmt$, we derive a more general lemma.
It illustrates that backtracking along any suitable direction is well-defined.
\begin{Lemma}
  Suppose 
  \cref{assumption:feasible_set_compact_convex,assumption:f_m_twice_continuously_differentiable}
  hold.
  For $\xitat \in \mathbb{R}^n$, let $\ve d$ be a descent direction for $\vemt$ 
  and let $q\in \{1,\ldots,k\}$ be any objective index and $\bar \sigma > 0$.
  Then there is an integer $j\in \mathbb{N}_0$ such that
  \begin{equation}
  \Psi
    \left(
    \xitat + \frac{\backtrackconst^j \bar \sigma}{\norm{\ve d}} \ve d
    \right)
    \le
    \Psi( \xitat )
    - \frac{\armijofactor \backtrackconst^j \bar \sigma}{\norm{\ve d}}
    w
    ,\qquad
    \armijofactor, \backtrackconst\in (0,1),
  \label{eqn:general_backtracking_condition}
  \end{equation}
  where, again, we have used the shorthand notation
  $
    w = -\max_{\ell=1,\ldots,k}\left\langle
    \gradmt_\ell(\xitat), \ve d
    \right\rangle > 0
  $
  and $\Psi$ is either some specific model, $\Psi = m_\ell$, 
  or the maximum value, $\Psi = \Phimt$.
  % \spremark{
  %    That is, the descent is guaranteed either for a specific component or for 
  %    the objective with the maximal value.
  %}
  Moreover, if we define the step 
  % $\stept$ as
  $
  \stept = \frac{\backtrackconst^j \bar \sigma}{\norm{\ve d}} \ve d
  $
  for the \textbf{smallest} $j\in \mathbb{N}_0$ satisfying 
  \eqref{eqn:general_backtracking_condition}, 
  then there is a constant $\sufficientconstantm\in (0,1)$ such that
  \begin{equation}
      \Psi(\xitat ) - \Psi\left(
          \xitat + \stept
      \right)
      \ge
          \sufficientconstantm w
          \min
          \left\{
          \frac{
            w
          }{
            \norm{\ve d}^2\normconst\hessboundmt
          },
          \frac{\bar \sigma}{\norm{\ve d}}
          \right\}.
      \label{eqn:general_backtracking_sufficient_decrease}
  \end{equation}
  \label{theorem:backtrack_one_function}
\end{Lemma}

\begin{proof}
The first part can be derived from the fact that $\ve d$ is a descent direction, 
see e.g. \cite{fukuda_drummond}.
However, we will use the approach from \cite{villacorta} to 
also derive the bound \eqref{eqn:general_backtracking_sufficient_decrease}.
With Taylor's Theorem we obtain
\begin{align}
    &\Psi
        \left(
            \xitat + \frac{\backtrackconst^j \bar \sigma}{\norm{\ve d}} \ve d
        \right)
    =
%\max_{q=1,\ldots,\ell}
        m_\ell
        \left(
            \xitat + \frac{\backtrackconst^j \bar \sigma}{\norm{\ve d}} \ve d
        \right)
            \qquad \text{(for some $\ell\in\{1,\ldots,k\})$}
    \nonumber\\&=
        %\max_{q = 1,\ldots,k} 
        %\left\{
            \mitat_\ell(\xitat)
            + \frac{\backtrackconst^j\bar{\sigma}}{\norm{{\ve d}}} 
                \langle \gradmt_\ell(\xitat), {\ve d} \rangle
            + \frac{(\backtrackconst^j\bar{\sigma})^2 }{2\norm{{\ve d}}^2} 
                \langle {\ve d}, \hessianmt{\ell}(\ve{\upxi}_\ell) {\ve d} \rangle
        %\right\}
    \nonumber\\&\le
        \Psi(\xitat)
        + \max_{q=1,\ldots,k} \frac{\backtrackconst^j\bar{\sigma}}{\norm{{\ve d}}} \langle \gradmt_q(\xitat), {\ve d} \rangle
        + \max_{q=1,\ldots,k} \frac{(\backtrackconst^j\bar{\sigma})^2 }{2\norm{{\ve d}}^2} \langle {\ve d}, \hessianmt{q}(\ve{\upxi}_q) {\ve d} \rangle
    \nonumber\\&\stackrel{
        \eqref{eqn:descent_direction_problemm},
        %\text{C.S.}, 
        \eqref{eqn:hessian_upper_bounds_ita}
    }\le
    \Psi(\xitat)
        - \frac{\backtrackconst^j\bar{\sigma}}{\norm{{\ve d}}} w
        + \frac{(\backtrackconst^j\bar{\sigma})^2 }{2} \normconst \hessboundmt.
%\label{eqn:modified_cauchy_step_first_estimate}
\label{eqn:general_backtracking_taylor_upper_estimate}
\end{align}
In the last line, we have additionally used the Cauchy-Schwarz inequality.\\
For a constructive proof, suppose now that \eqref{eqn:general_backtracking_condition}
 is violated for some $j\in \mathbb{N}_0$, i.e.,
$$
\Psi
	\left(
	\xitat + \frac{\backtrackconst^j \bar \sigma}{\norm{\ve d}} \ve d
	\right)
	>
	\Psi( \xitat )
	- \frac{\armijofactor \backtrackconst^j \bar \sigma}{\norm{\ve d}}
    w
    .
$$
Plugging in \eqref{eqn:general_backtracking_taylor_upper_estimate} for the LHS 
and substracting $\Psi(\xitat)$ then leads to
% \begin{equation*}
% 	- \frac{\backtrackconst^j\bar{\sigma}}{\norm{{\ve d}}} w
% 	+ \frac{(\backtrackconst^j\bar{\sigma})^2 }{2} \normconst \hessboundmt
% 	>
% 	- \frac{\armijofactor \backtrackconst^j \bar \sigma}{\norm{\ve d}}	w
% \end{equation*} 
%   % \Rightarrow \quad&&
% 	% - \frac{1}{\norm{\ve d}}w + \frac{\backtrackconst^j\bar{\sigma}}{2} \normconst \hessboundmt
% 	% & > - \frac{\armijofactor}{\norm{\ve d}} w
%   % \nonumber
% and thus
\begin{equation*}
	\backtrackconst^j
	>
	\frac{
		2(1-\armijofactor) w
		}{
			\norm{\ve d}\bar{\sigma}\normconst\hessboundmt
		} ,
\end{equation*}
where the right hand side is positive and completely independent of $j$.
Since $\backtrackconst\in(0,1)$, there must be a $j^*\in \mathbb{N}_0, j^*>j,$ for which
$
\backtrackconst^{j^*} \le
\dfrac{
	2(1-\armijofactor) w
	}{
		\norm{\ve d}\bar{\sigma}\normconst\hessboundmt
	}
$
so that \eqref{eqn:general_backtracking_condition} must also be fulfilled
for this $\backtrackconst^{j^*}$.

Analogous to the proof of \cite[Lemma 4.2]{villacorta}
we can now derive the constant $\sufficientconstantm$ from
\eqref{eqn:general_backtracking_sufficient_decrease} as 
$
  \sufficientconstantm = \min \{ 2\backtrackconst(1-\armijofactor), \armijofactor \}.
$

% Let $j\in \mathbb{N}_0$ be the smallest integer satisfying \eqref{eqn:general_backtracking_condition}.\\
% If $j\ge 1$, then $j-1$ violates the condition and we find
% \begin{equation*}
%     \backtrackconst^{j} 
%     = \backtrackconst\backtrackconst^{j-1}
%     >
%         \backtrackconst
%         \frac{
%             2(1-\armijofactor) w
%             }{
%                 \norm{\ve d}\bar{\sigma}\normconst\hessboundmt
%             }.
%     %\label{eqn:backtrackconst_pow_j_strict_lower_bound}
% \end{equation*}
% If we plug this into the backtracking condition \eqref{eqn:general_backtracking_condition} then
% $$
% \Psi(\xitat ) - \Psi\left(
% \xitat + \stept
% \right)
% \ge
% \frac{
% 	2\backtrackconst(1-\armijofactor) w^2
% 	}{
% 		\norm{\ve d}^2\normconst\hessboundmt
% 	}.
% $$
% If otherwise $j=0$, then the RHS in \eqref{eqn:general_backtracking_condition} reduces to
% $
% \dfrac{\armijofactor\bar{\sigma}w}{\norm{\ve d}}
% $ so that finally,
% \begin{align*}
%   \Psi(\xitat ) - \Psi\left(
%   \xitat + \stept
%   \right)
%   &\ge
%   w
%   \min \left\{
%   \frac{
%     2\backtrackconst(1-\armijofactor) w
%     }{
%       \norm{\ve d}^2\normconst\hessboundmt
%     },
%     \frac{\armijofactor\bar{\sigma}}{\norm{\ve d}}
%   \right\}
%   \\
%   &\ge
%   \sufficientconstantm
%   w
%   \min\left\{
%     \frac{ w}{\norm{\ve d}^2\normconst\hessboundmt},
%     \frac{\bar \sigma}{\norm{\ve d}}
%   \right\},
% \end{align*}
% where we have defined
% $
%   \sufficientconstantm = \min \{ 2\backtrackconst(1-\armijofactor), \armijofactor \} .
%   %\label{eqn:def_sufficient_constant}
% $.
\end{proof}

\cref{theorem:backtrack_one_function} applies naturally to the step along $\dmt$:
\begin{Definition}
    For $\xitat \in \ballt$ let $\dmt$ be a solution to $\eqref{eqn:descent_direction_problemm}$ 
    and define the \emph{modified Pareto-Cauchy step} as
    \begin{equation*}
        \modscauchyt := \backtrackconst^j \bar{\sigma} \frac{\dmt}{\normdmt},
  \end{equation*}
    where
    again $\bar \sigma$ as in \eqref{eqn:bar_sigma_def} and
    $j\in \mathbb{N}_0$ is the smallest integer that satisfies
    \begin{equation}
        \Phimt( \xitat + \modscauchyt ) \le \Phimt(\xitat) -
        \frac{
            \armijofactor \backtrackconst^j \bar{\sigma}
        }{
            \normdmt
        }
        \omegamt{ \xitat }
        \label{eqn:backtracking_condition}
    \end{equation}
    for predefined constants $\armijofactor, \backtrackconst \in (0,1)$.
    \label{mydef:mod_pareto_cauchy}
\end{Definition}

The definition of $\bar \sigma$ ensures, that $\xitat + \modscauchyt$ is contained 
in the current trust region $\ballt$.
Furthermore, these steps provide a sufficient decrease 
very similar to \eqref{eqn:sufficient_decrease_true_cauchy_point}:

\begin{Corollary}
    Suppose \cref{assumption:feasible_set_compact_convex,assumption:f_m_twice_continuously_differentiable} hold.
    For the step $\modscauchyt$ the following statements are true:
    \begin{enumerate}
    \item A $j\in \mathbb{N}_0$ as in \eqref{eqn:backtracking_condition} exists.
    \item There is a constant $\sufficientconstantm \in (0,1)$ such that the modified 
    Pareto-Cauchy step $\modscauchyt$ satisfies
        \begin{equation*}
            \Phimt(\xitat) - \Phimt( \xitat + \modscauchyt )
            \ge
            \sufficientconstantm \omegamt{  \xitat  } \min
            \left\{
                \frac{\omegamt{ \xitat }}{\normconst\hessboundmt}, \radiust, 1
            \right\}.
      \end{equation*}
    \end{enumerate}
    \label{theorem:sufficient_decrease_mod_cauchy_point}
\end{Corollary}

\begin{proof}
If $\xitat$ is critical, then the bound is trivial.
Otherwise, the existence of a $j$ satisfying \eqref{eqn:backtracking_condition} follows 
from \cref{theorem:backtrack_one_function} for $\Psi = \Phimt$.
The lower bound on the decrease follows immediately from 
$\bar \sigma \ge \min\left\{ \normdmt, \radiust \right\}$.
\end{proof}

%\subsubsection*{Strict Modified Pareto-Cauchy Step}
From \cref{theorem:backtrack_one_function} it follows that 
the backtracking condition \eqref{eqn:backtracking_condition} can be 
modified to explicitly require a decrease in \emph{every} objective:
\begin{Definition}
	Let $j\in \mathbb{N}_0$ the smallest integer satisfying
	\begin{equation*}
		\min_{\ell=1,\ldots,k}
		\left\{
	  	\mitat_\ell(\xitat) - \mitat_\ell\left(\xitat + \backtrackconst^j\bar{\sigma}\frac{\dmt}{\normdmt}\right)
		\right\}
		\ge
		\frac{\armijofactor\backtrackconst^j\bar{\sigma}}{\normdmt}\omegamt{ \xitat }.
		%\label{eqn:strict_mod_cauchy}
  \end{equation*}
  We define the \emph{strict} modified Pareto-Cauchy point as 
  $\strictx = \xitat + \stricts $ 
  and the corresponding step as
  $
  \stricts = \backtrackconst^j\bar{\sigma}\dfrac{\dmt}{\normdmt}
  $.
	% \begin{equation*}
	% 	\strictx = \xitat + \stricts , \quad \stricts = \backtrackconst^j\bar{\sigma}\frac{\dmt}{\normdmt}.
	% 	%\label{eqn:def_strictx}
  % \end{equation*}
	\label{mydef:strict_pareto_cauchy}
\end{Definition}

\begin{Corollary}
	Suppose \cref{assumption:feasible_set_compact_convex,assumption:f_m_twice_continuously_differentiable} hold.
	\begin{enumerate}
	\item The strict modified Pareto-Cauchy point exists, the backtracking is finite.
	\item There is a constant $\sufficientconstantm \in (0,1)$ such that
	\begin{equation}
    \min_{\ell=1,\ldots,k}
    \left\{
      \mitat_\ell(\xitat) - \mitat_\ell\left( \strictx \right)
    \right\}
    \ge
    \sufficientconstantm
    \omegamt{  \xitat  } \min
    \left\{
      \frac{\omegamt{ \xitat }}{\normconst\hessboundmt}, \radiust, 1
    \right\}.
    \label{eqn:strict_sufficient_decrease}
  \end{equation}
  \end{enumerate}
  \label{theorem:strict_sufficient_decrease}
\end{Corollary}

\begin{Remark}
  In the preceding subsections, we have shown descent steps along the model steepest descent direction.
  Similar to the single objective case we do not necessarily have to use the steepest descent direction
  and different step calculation methods are viable.
  For instance, \citet*{thomann} use the well-known Pascoletti-Serafini 
  scalarization to solve the subproblem \eqref{eqn:mopmt}.
  We refer to their work and \cref{appendix:section_steps} to see how this method can be related 
  to the steepest descent direction.
\end{Remark}
%\textcolor{magenta}{In the Appendix we also briefly discuss a nonlinear conjugate gradient method.}

%%% (start) input of file  /home/manuelbb/Desktop/latex_cleanup/trm_combined/sections/4_steps/sufficient_decrease.tex  %%%%
\subsection{Sufficient Decrease for the Original Problem}

In the previous subsections, we have shown how to compute steps $\stept$ to achieve
a sufficient decrease in terms of $\Phimt$ and $\omegamt{\bullet}$.
For a descent step $\stept$ the bound is of the form
\begin{equation}
  \Phimt( \xitat ) - \Phimt( \xitat + \stept )
  \ge
  \sufficientconstantm \omegamt{ \xitat }
  \min
  \left\{
    \frac{ \omegamt{\xitat} }{ \normconst \hessboundmt },
    \radiust,
    1
  \right\}, \quad
  \sufficientconstantm \in (0,1),
  \label{eqn:sufficient_decrease_prototype}
\end{equation}
and thereby very similar to the bounds for the scalar projected gradient trust region method \cite{conn_trust_region_2000}.
By introducing a slightly modified version of $\omegamt{\bullet}$, 
we can transform \eqref{eqn:sufficient_decrease_prototype}
into the bound used in \cite{thomann} and \cite{villacorta}.

\begin{Lemma}
  If $\pi(t,\xitat)$ is a criticality measure for some multiobjective problem, then
  $
    \tilde \pi(t, \xitat) = \min \left\{1, \pi(t, \xitat) \right\}
  $
  is also a criticality measure for the same problem.
  \label{theorem:modified_criticality_measure}
\end{Lemma}

\begin{proof}
  We have
  $
  0 \le \tilde \pi( t, \xitat) \le \pi(t, \xitat)$.
  Thus, $\tilde \pi \to 0$ whenever $\pi\to 0$.
  The minimum of uniformly continuous functions is again uniformly continuous.
\end{proof}

We next make another standard assumption on the class of surrogate models.
% unconstrained assumption
\begin{Assumption}
  The norm of all model hessians is uniformly bounded above on $\feas$, 
  i.e., there is a positive constant $\hessboundm$ such that
  \begin{equation*}
    \left\| \hessianmt{\ell} ( \ve x) \right\|_F \le \hessboundm \qquad 
    \forall \ell = 1, \ldots, k, \forall \ve x \in \ballt, \; 
    \forall t \in \mathbb{N}_0.
  \end{equation*}
  W.l.o.g., we assume
  \begin{equation}
    \hessboundm \cdot \normconst > 1, \quad \text{with $\normconst$ as in \eqref{eqn:normconst}.}
	  \label{eqn:hessbound_times_normconst}
  \end{equation}
  \label{assumption:model_hessian_bounded_above}
\end{Assumption}

\begin{Remark}
  From this assumption it follows that the model gradients are then Lipschitz as well.
  Together with \cref{theorem:omega_uniformly_continuous}, 
  we then know that $\omegamt{\bullet}$ is a criticality measure for \eqref{eqn:mopmt}.
\end{Remark}

Motivated by the previous remark, we will from now on refer to the following functions
\begin{equation}
\begin{aligned}
	\modomega{\ve x}
		:=
			\min \left\{ \omga{ \ve x }, 1 \right\}
		\, \text{ and } \,
	\modomegamt{ \ve x }
		:=
		\min \left\{ \omegamt{ \ve x }, 1 \right\} &\forall t = 0, 1,\ldots
\end{aligned}
\label{eqn:modomega_definition}
\end{equation}
We can thereby derive the sufficient decrease condition in “standard form”:

\begin{Corollary}
  Under \cref{assumption:model_hessian_bounded_above}, suppose that for $\xitat$ 
  and some descent step $\stept$ the bound \eqref{eqn:sufficient_decrease_prototype} holds.
  For the criticality measure $\modomegamt{\bullet}$ it follows that
  \begin{equation}
  \Phimt( \xitat ) - \Phimt( \xitat + \stept )
  \ge
  \sufficientconstantm
  \modomegamt{ \xitat }
  \min
  \left\{
    \frac{ \modomegamt{\xitat} }{ \normconst \hessboundm },
    \radiust
  \right\}.
  \label{eqn:sufficient_decrease_modified_prototype}
  \end{equation}
  \label{theorem:sufficient_decrease_modified_prototype}
\end{Corollary}

\begin{proof}
$\modomegamt{\bullet}$ is a criticality measure due to
\cref{assumption:model_hessian_bounded_above,theorem:modified_criticality_measure}.
Further, from \eqref{eqn:modomega_definition} and \eqref{eqn:hessbound_times_normconst} it follows that
$$
  \frac{ \modomegamt{\xitat} }{ \normconst \hessboundm }
  \le
    \frac{ 1 }{ \normconst \hessboundm }
  \le
    1
$$
and if we plug this into \eqref{eqn:sufficient_decrease_prototype} 
we obtain \eqref{eqn:sufficient_decrease_modified_prototype}.
\end{proof}

To relate the RHS of \eqref{eqn:sufficient_decrease_modified_prototype} to 
the criticality $\omga{\bullet}$ of the original problem, we require another assumption.
\begin{Assumption}
  There is a constant $\omegaconstant > 0$ such that
  \begin{equation*}
    \left|
      \omegamt{ \xitat } - \omga{  \xitat  }
    \right|
    \le
    \omegaconstant \omegamt{  \xitat  }.
    %\label{eqn:omegaconstant_definition}
  \end{equation*}
  \label{assumption:omegaconstant}
\end{Assumption}
This assumption is also made by \citet{thomann} and can easily be 
justified by using fully linear surrogate models and a bounded trust 
region radius in combination with the a criticality test, 
see \cref{theorem:omegaconst_criticality}.\\
\cref{assumption:omegaconstant} can be used to formulate the next two lemmata relating
the model criticality and the true criticality.
They are proven in \cref{appendix:modomega}. 
From these lemmata and \cref{theorem:sufficient_decrease_modified_prototype} 
the final result, \cref{theorem:sufficient_decrease_prototype_modomega},
easily follows.

\begin{Lemma}
  If \cref{assumption:omegaconstant} holds, then it holds for $\modomegamt{ \bullet }$ and $\modomega{ \bullet }$
  from \eqref{eqn:modomega_definition} 
  %\spout{too }
  that
  \begin{equation*}
    \left|
      \modomegamt{ \xitat } - \modomega{  \xitat  }
    \right|
    \le
    \omegaconstant \modomegamt{  \xitat  }.
  \end{equation*}
  \label{theorem:omegaconst_mod}
\end{Lemma}

% \begin{proof}
%   The proof can be found in \cref{appendix:modomega}
% \end{proof}

\begin{Lemma}
  From \cref{assumption:omegaconstant} it follows that
  $$
	 \modomegamt{  \xitat  } \ge \frac{1}{\omegaconstant + 1} \modomega{  \xitat  }
    \quad \text{with $(\omegaconstant+1)^{-1}\in (0,1).$}
  $$
\label{theorem:omegaconst}
\end{Lemma}

% \begin{proof}
%   The proof can again be found in \cref{appendix:modomega}.
% \end{proof}

\begin{Corollary}
  Suppose that \cref{assumption:model_hessian_bounded_above,assumption:omegaconstant} 
  hold and that $\xitat$ and $\stept$ satisfy \eqref{eqn:sufficient_decrease_modified_prototype}.
  Then
  \begin{align}
    \Phimt( \xitat ) - \Phimt( \xitat + \stept )
    &\ge
      \sufficientconstant
      \modomega{ \xitat }
      \min
      \left\{
        \frac{ \modomega{\xitat} }{ \normconst \hessboundm },
        \radiust
      \right\},     
    \label{eqn:sufficient_decrease_prototype_modomega}
  \end{align}
  where  
  $
    \sufficientconstant = \frac\sufficientconstantm{1 + \omegaconstant}
    \in (0,1)
  $.
  \label{theorem:sufficient_decrease_prototype_modomega}
\end{Corollary}

%%% (start) input of file  /home/manuelbb/Desktop/latex_cleanup/trm_combined/sections/6_convergence/assumptions.tex  %%%%

\section{Convergence}
\label{section:convergence}

\subsection{Preliminary Assumptions and Definitions}

\noindent To prove convergence of \cref{algo:algorithm1} we first have to 
make sure that at least one of the objectives is bounded from below:
\begin{Assumption}
	The maximum $\max_{\ell=1,\ldots,k} f_\ell(\ve x)$ of all objective functions 
	is bounded from below on $\feas$.
	\label{assumption:objectives_bounded}
\end{Assumption}

\noindent 
To be able to use $\modomega{\bullet}$ as a criticality measure and to 
refer to fully linear models, we further require:
\begin{Assumption}
	The objective $\ve f\colon \mathbb{R}^n\to \mathbb{R}^k$ is continuously differentiable in an open domain containing $\feas$ 
	and has a Lipschitz continuous gradient on $\feas$.
	\label{assumption:objective_gradients_lipschitz}
\end{Assumption}

\noindent We summarize the assumptions on the surrogates as follows:
\begin{Assumption}
	The surrogate model functions $m_1\itat,\ldots,m_k\itat$ belong 
	to a fully linear class $\mathcal M$ as defined in \cref{mydef:fully_linear}.
%	For \cref{algo:algorithm0} we actually require all models to be fully linear.
	For each objective index $\ell\in\{1,\ldots,k\}$, 
	the error constants are then denoted by 
	$\fullylinearconstant_\ell$ and $\fullylinearconstantdf_\ell$.
%	We further assume every $m_\ell\itat$ to be twice continuously differentiable and the hessians to be uniformly bounded by $\hessboundm$ as in \eqref{eqn:hessbound_times_normconst}.
\label{assumption:surrogates}
\end{Assumption}

\noindent For the subsequent analysis we define component-wise maximum constants as
\begin{equation}
	\fullylinearconstant := \max_{\ell = 1, \ldots, k} \fullylinearconstant_\ell, \quad
	\fullylinearconstantdf := \max_{\ell=1,\ldots,k} \fullylinearconstantdf_\ell.
	\label{eqn:fullylinearconstants_max}
\end{equation}

\noindent We also wish for the descent steps to fulfill a sufficient decrease condition for
the surrogate criticality measure as discussed in \cref{section:steps}.
\begin{Assumption}
	For all $t\in \mathbb{N}_0$ the descent steps $\stept$ are assumed to fulfill both $\xitat + \stept \in \ballt$ and \eqref{eqn:sufficient_decrease_modified_prototype}.
	\label{assumption:sufficient_decrease}
\end{Assumption}

\noindent Finally, to avoid a cluttered notation when dealing with 
subsequences we define the following shorthand notations:
\begin{equation*}
	\modshortmit{t} := \modomegamit{ t } ,
	\;
	\modshortit{t} := \modomegait{ t }\quad
	\forall t\in \mathbb{N}_0.
\end{equation*}

\subsection{Convergence of \cref{algo:algorithm1}}

In the following we prove convergence of \cref{algo:algorithm1} to Pareto critical points. 
We account for the case that no criticality test is used, i.e., $\epscrit = 0$.
We then require all surrogates to be fully linear in each iteration
and need \cref{assumption:omegaconstant}.
The proof is an adapted version of the scalar case in \cite{conn:trm_framework}.\\
It is also similar to the proofs for the multiobjective algorithms in \cite{thomann,villacorta}.
However, in both cases, no criticality test is employed, 
there is no distinction between successful and acceptable iterations ($\nuaccept = \nusuccess$) 
and interpolation at $\xitat$ by the surrogates is required.
We indicate notable differences when appropriate.

We start with two results concerning the criticality test in \cref{algo:algorithm1}.

\begin{Lemma}
	Outside the \criticalityRoutine, \cref{assumption:omegaconstant} is fulfilled 
	if the model $\vemt$ is fully-linear 
	(and if $\radiust \le \radiusmax < \infty$).
	\label{theorem:omegaconst_criticality}
\end{Lemma}

\begin{proof}
Let $\ell, q\in \{1,\ldots,k\}$ and $\ve d_\ell, \ve d_q\in \feas -\xitat$ be 
solutions of \eqref{eqn:descent_direction_problem1} and \eqref{eqn:descent_direction_problemm}
respectively such that
\begin{align*}
  \omegamt{\xitat} = -\langle \gradmt_\ell(\xitat), \ve d_\ell \rangle,\;
  \omga{\xitat} = -\langle \gradf_q(\xitat), \ve d_q \rangle.
\end{align*}
If $\omegamt{\xitat}\ge \omga{\xitat}$, then, 
using Cauchy-Schwarz and $\norm{\ve d_\ell}\le 1$,
\begin{align*}
  \left|
    \omegamt{\xitat} - \omga{\xitat}
  \right|
  &=
    \langle \gradf_q(\xitat), \ve d_q \rangle
    -\langle \gradmt_\ell(\xitat), \ve d_\ell \rangle
  \\
  &\stackrel{\text{df.}}\le
    \langle \gradf_q(\xitat), \ve d_\ell \rangle
    -\langle \gradmt_q(\xitat), \ve d_\ell \rangle
  \\
  &\le
  \norm{
    \gradf_q(\xitat) - \gradmt_q(\xitat)
  }_2,
\end{align*}
and if $\omegamt{\xitat}< \omga{\xitat}$, we obtain
\begin{align*}
  \left|
    \omegamt{\xitat} - \omga{\xitat}
  \right|
  &\le
  \norm{
    \gradmt_\ell(\xitat) - \gradf_\ell(\xitat)
  }_2.
\end{align*}
Because $\vemt$ is fully linear, it follows that
\begin{align*}
  \left|
    \omegamt{\xitat} - \omga{\xitat}
  \right|
  &\le
  \sqrt{\normconst}\fullylinearconstantdf \radiust,
  \qquad\text{with $\fullylinearconstantdf$ from \eqref{eqn:fullylinearconstants_max}.}
\end{align*}

If we just left \criticalityRoutine, then the model is fully linear 
for $\radiust$ due to \cref{theorem:fully_linear_in_larger_region} and we have 
$\radiust \le \mu\modomegamt{\xitat} \le \mu \omegamt{\xitat}$.
If we otherwise did not enter \criticalityRoutine in the first place,
it must hold that $\omegamt{\xitat} \ge \epscrit$ and
$$
  \radiust \le \radiusmax = \frac{\radiusmax}{\epscrit}\epscrit
  \le
    \frac{\radiusmax}{\epscrit}\omegamt{\xitat}
$$
and thus
$$
  \left|
    \omegamt{\xitat} - \omga{\xitat}
  \right|
  \le
    \omegaconstant \omegamt{\xitat},
    \quad
    \omegaconstant = \sqrt{\normconst}
    \fullylinearconstantdf
    \max\left\{
      \mu, \epscrit^{-1}\radiusmax
    \right\}
  > 0.
$$
\end{proof}

In the subsequent analysis, we require mainly steps with fully linear models to 
achieve sufficient decrease for the true problem.
Due to \cref{theorem:omegaconst_criticality}, 
we can dispose of \cref{assumption:omegaconstant} by using the criticality routine:

\begin{Assumption}
	Either $\epscrit > 0$ or \cref{assumption:omegaconstant} holds.
	\label{assumption:omegaconstant_relaxed}
\end{Assumption}

We have also implicitly shown the following property of the criticality measures.
\begin{Corollary}
  If $\vemt$ is fully linear for $\ve f$ with $\fullylinearconstantdf > 0$ as in \eqref{eqn:fullylinearconstants_max} then
  \begin{align*}
    \left|
      \modomegamt{\xitat} - \modomega{\xitat}
    \right|
    \le
    \left|
      \omegamt{\xitat} - \omga{\xitat}
    \right|
    \le \sqrt{c} \fullylinearconstantdf \radiust.
  \end{align*}
  \label{theorem:omegas_fully_linear}
\end{Corollary}

\begin{Lemma}
  If $\xitat$ is not critical for the true problem \eqref{eqn:mop}, 
  i.e. $\modomega{\xitat} \ne 0$, 
  then \allowbreak\criticalityRoutine
  will terminate after a finite number of iterations.
	\label{theorem:inifinite_criticality_loop}
\end{Lemma}

\begin{proof}
%The proof is inspired by a similar result 
%for single objective optimization \cite[Lemma 5.1]{conn:trm_framework}.
At the start of \criticalityRoutine, we know that $\vemt$ is 
not fully linear or $\radiust > \mu\modomegamt{\xitat}$.
For clarity, we denote the first model by $\vemt_0$ and define $\Delta_0 = \radiust$.
We then ensure that the model is made fully linear on $\radiust_1 = \Delta_0$ 
and denote this fully linear model by $\vemt_1$.
If afterwards $\radiust_1 \le \mu\modomegamtcrit{1}{\xitat}$, 
then \criticalityRoutine terminates. \\
Otherwise, the process is repeated: 
the radius is multiplied by $\alpha\in (0,1)$ so that in the $j$-th iteration
we have  $\radiust_j = \alpha^{j-1} \Delta_0$ and 
$\vemt_j$ is made fully linear on $\radiust_j$ until
\begin{equation*}
  \radiust_j =\alpha^{j-1} \Delta_0 \le \mu \modomegamtcrit{j}{\xitat}.
  %\label{eqn:criticality_end}
\end{equation*}
The only way for \criticalityRoutine to loop infinitely is
\begin{align}
  \modomegamtcrit{j}{\xitat} &< \frac{\alpha^{j-1} \Delta_0}{\mu} \qquad \forall j\in \mathbb{N}.
  \label{eqn:inifinity_condition}
\end{align}
Because $\vemt_j$ is fully linear on $\alpha^{j-1} \Delta_0$, we know from \cref{theorem:omegas_fully_linear} that
$$
  \left| \modomegamtcrit{j}{\xitat} - \modomega{\xitat} \right|
  \le
  %\sqrt{\normconst}\fullylinearconstantdf\radiust_j =
  \sqrt{\normconst}\fullylinearconstantdf \alpha^{j-1} \Delta_0 \qquad \forall j\in \mathbb{N}.
$$
Using the triangle inequality together with \eqref{eqn:inifinity_condition} gives us
$$
  \left|
    \modomega{\xitat}
  \right|
  \le
    \left|
      \modomegamtcrit{j}{\xitat} - \modomega{\xitat}
    \right|
    +
    \left|
      \modomegamtcrit{j}{\xitat}
    \right|
  \le
    \left( \mu^{-1} + \sqrt{\normconst} \fullylinearconstant \right)
    \alpha^{j-1} \Delta_0 \quad \forall j\in \mathbb{N}.
$$
As $\alpha\in (0,1)$, this implies $\modomega{\xitat} = 0$ and $\xitat$ is hence critical.
\end{proof}

We next state another auxiliary lemma that we need for the convergence proof.

\begin{Lemma}
	Suppose \cref{assumption:objective_gradients_lipschitz,assumption:surrogates} hold.
	For the iterate $\xitat$ let $\stept\in \mathbb{R}^n$ be a any step with
	$\xtrialt = \xitat + \stept \in \ballt$.
	If $\vemt$ is fully linear on $\ballt$ then it holds that
	\begin{equation*}
		\left|
			\Phif(\xtrialt) - \Phimt(\xtrialt)
		\right|
		\le
		\fullylinearconstant \left( \radiust\right)^2.
	\end{equation*}
	\label{lemma:absvalue_Phi_bounded}
\end{Lemma}

\begin{proof}
	The proof follows from the definition of $\Phif$ and $\Phimt$ and the full linearity of $\vemt$.
	It can be found in \cite[Lemma 4.16]{thomann}.
\end{proof}

Convergence of \cref{algo:algorithm1} is proven by showing that in certain situations, 
the iteration must be acceptable or successful as defined in \cref{mydef:iter_categories}.
This is done indirectly and relies on the next two lemmata.
They use the preceding result to show 
that in a (hypothetical) situation where no Pareto-critical point is approached, 
the trust region radius must be bounded from below.

\begin{Lemma}
	Suppose 
	\cref{assumption:feasible_set_compact_convex,assumption:model_hessian_bounded_above,assumption:objective_gradients_lipschitz,assumption:surrogates,assumption:sufficient_decrease,assumption:feasible_set_compact_convex} 
	hold.
	If $\xitat$ is not Pareto-critical for \eqref{eqn:mopmt} and $\vemt$ is fully linear on $\ballt$ and
	\begin{equation*}
		\radiust \le
		\frac{
			\sufficientconstantm ( 1 - \nusuccess) \modomegamt{ \xitat }
		}{
			2\lambda
		},
		\quad \text{where $\lambda = \max\left\{ \fullylinearconstant, \normconst\hessboundm \right\}$ and 
		$\sufficientconstantm$ as in \eqref{eqn:sufficient_decrease_modified_prototype},}
		%\label{eqn:upper_bound_on_radiust}
	\end{equation*}
	then the iteration is successful, that is, $t\in \successindices$ and $\Delta^{t+1}\ge \radiust$.
	\label{theorem:small_radius_very_successful}
\end{Lemma}

\begin{proof}
The proof is very similar to \cite[Lemma 5.3]{conn:trm_framework} 
and \cite[Lemma 4.17]{thomann}. 
In contrast to the latter, we use the surrogate problem and do not
require interpolation at $\xitat$:

By definition we have $\sufficientconstantm(1-\nusuccess) < 1$ and hence it follows from
\cref{assumption:sufficient_decrease,assumption:omegaconstant,theorem:sufficient_decrease_modified_prototype} that
\begin{align}
	\radiust
	&\le
		\frac{
			\sufficientconstantm ( 1 - \nusuccess) \modomegamt{ \xitat }
		}{
			2\lambda
		}
	\label{eqn:upper_bound_on_radiust_by_sufficient_constant}
	\\
	&\le
		\frac{\modshortmit{t}}{2\lambda}
	\le
		\frac{\modshortmit{t}}{2\normconst\hessboundm}
	\le
		\frac{\modshortmit{t}}{\normconst\hessboundm}.
		\nonumber
\end{align}
With \cref{assumption:sufficient_decrease} we can plug this into \eqref{eqn:sufficient_decrease_modified_prototype} and obtain
\begin{equation}
	\Phimt(\xitat) - \Phimt( \xtrialt )
	\ge
		{\sufficientconstantm} \modshortmit{t}
		\min \left\{
			\frac{
				\modshortmit{t}
			}{
				\normconst\hessboundm
			},
			\radiust
		\right\}
	\ge
		{\sufficientconstantm} \modshortmit{t}\radiust.
		\label{eqn:sufficient_decrease_radius_only}
\end{equation}
Due to \cref{assumption:surrogates} we can take the definition \eqref{eqn:rho_definition} and estimate
\begin{align*}
	\left|
		\rho\ita{t} - 1
	\right|
	&=
	\left|
		\frac{
			\Phif(\xitat) - \Phif(\xtrialt) - (\Phimt( \xitat ) - \Phimt( \xtrialt )
		}{
			\Phimt(\xitat) - \Phimt( \xtrialt)
		}
	\right|
	\\
	&\le
		\frac{
			\left| \Phif(\xitat)  - \Phimt( \xitat ) \right|
			+
			\left|  \Phimt( \xtrialt ) - \Phif(\xtrialt) \right|
		}{
			\left|\Phimt(\xitat) - \Phimt( \xtrialt)\right|
		}
	\\&\stackrel{
		\text{\cref{lemma:absvalue_Phi_bounded},
		\eqref{eqn:sufficient_decrease_radius_only}}
	}\le
	\frac{
		2 \fullylinearconstant \left(\radiust\right)^2
	}{
		\sufficientconstantm\modshortmit{t}\radiust
	}
	\le
	\frac{
		2 \lambda \radiust 
	}{
	\sufficientconstantm\modshortmit{t}
	}
	\stackrel{\eqref{eqn:upper_bound_on_radiust_by_sufficient_constant}}\le
	1 - \nusuccess.
\end{align*}
Therefore $\rho\itat \ge \nusuccess$ and the iteration $t$ using step $\stept$ is successful.
\end{proof}

The same statement can be made for the true problem and $\modomega{\bullet}$:
\begin{Corollary}
	Suppose \cref{assumption:feasible_set_compact_convex,assumption:model_hessian_bounded_above,assumption:objective_gradients_lipschitz,assumption:surrogates,assumption:sufficient_decrease,assumption:feasible_set_compact_convex,assumption:omegaconstant_relaxed} hold.
	If $\xitat$ is not Pareto-critical for \eqref{eqn:mop} and $\vemt$ is fully linear on $\ballt$ and
	\begin{equation*}
		\radiust \le
		\frac{
			\sufficientconstant ( 1 - \nusuccess) \modomega{ \xitat }
		}{
			2\lambda
		},
		\quad \text{where }\lambda = 
		\max\left\{ 
			\fullylinearconstant, \normconst\hessboundm 
		\right\}, 
		\sufficientconstantm
		\text{ as in \eqref{eqn:sufficient_decrease_prototype_modomega},}
		%\label{eqn:upper_bound_on_radiust}
	\end{equation*}
	then the iteration is successful, that is $t\in \successindices$ and $\Delta^{t+1}\ge \radiust$.
	\label{theorem:small_radius_very_successful_true_prolem}
\end{Corollary}
\begin{proof}
The proof works exactly the same as for \cref{theorem:small_radius_very_successful}.
But due to \cref{assumption:omegaconstant_relaxed} we can use \cref{theorem:omegaconst_criticality} and
employ the sufficient decrease condition \eqref{eqn:sufficient_decrease_prototype_modomega}
for $\modomega{\bullet}$ instead.
\end{proof}

As in \cite[Lemma 5.4]{conn:trm_framework} and \cite[Lemma 4.18]{thomann}, 
it is now easy to show that when no Pareto-critical point of
\eqref{eqn:mopmt} is approached the trust region radius must be bounded:
\begin{Lemma}
	Suppose \cref{assumption:feasible_set_compact_convex,assumption:model_hessian_bounded_above,assumption:objective_gradients_lipschitz,assumption:surrogates,assumption:sufficient_decrease} hold and
 	that there exists a constant $\lbomega > 0$ such that $\modomegamt{\xitat} \ge \lbomega$ for all $t$.
	Then there is a constant $\lbradius > 0$ with
	$$
		\radiust \ge \lbradius \quad \text{for all $t\in \mathbb{N}_0$}.
	$$
\label{lemma:lbomega_lbradius}
\end{Lemma}

\begin{proof}
We first investigate the criticality step and assume $\epscrit > \modshortmit{t} \ge \lbomega$.
After we finish the criticality loop, we get an radius $\radiust$ so that 
$\radiust \ge \min\{ \radiust_*, \beta\modshortmit{t}\}$ and therefore 
$\radiust \ge \min\{\beta\lbomega, \radiust_*\}$ for all $t$.

Outside the criticality step, we know from \cref{theorem:small_radius_very_successful}
that whenever $\radiust$ falls below
$$
	\tilde \Delta :=
	\frac{
		\sufficientconstantm ( 1 - \nusuccess ) \lbomega
	}
	{2\lambda},
$$
iteration $t$ must be either model-improving or successful and 
hence $\radiusnext \ge \radiust$ and the radius cannot decrease 
until $\radius\ita{k} > \tilde \Delta$ for some $k > t$.
Because $\gammasmallest\in (0,1)$ is the severest possible 
shrinking factor in \cref{algo:algorithm1}, we therefore know that
$\radiust$ can never be actively shrunken to a value below $\gammasmallest \tilde \Delta$.

Combining both bounds on $\radiust$ results in
$$
	\radiust \ge \lbradius := 
	\min\{ \beta\lbomega, \gammasmallest \tilde \Delta, \radius\ita0_* \} 
		\qquad \forall t\in \mathbb{N}_0,
$$
where we have again used the fact, that $\radiust_*$ cannot be reduced further 
if it is less than or equal to $\tilde \Delta$ due to the update mechanism in 
\cref{algo:algorithm1}.
\end{proof}

We can now state the first convergence result:
\begin{Theorem}
	Suppose that 
	\cref{assumption:model_hessian_bounded_above,assumption:objective_gradients_lipschitz,assumption:surrogates,assumption:sufficient_decrease,assumption:feasible_set_compact_convex} hold.
 	If \cref{algo:algorithm1} has only a finite number $0\le |\successindices| < \infty$ of successful iterations $\successindices = \{t\in \mathbb{N}_0: \rho\ita{t} \ge \nusuccess\}$ then
	$$
		\lim_{t\to \infty} \modomega{\xitat} = 0.
	$$
	\label{theorem:finitely_many_acceptable_iterations_convergence}
\end{Theorem}

\begin{proof}
%The proof is adapted from the single objective case \cite[Lemma 5.5]{conn:trm_framework} and 
%therefore differs a bit from the proof in \cite[Lemma 4.20]{thomann}.\\
If the criticality loop runs infinitely, then the result 
follows from \cref{theorem:inifinite_criticality_loop}.

Otherwise, let $t_0$ any index larger than the last successful index 
(or $t_0 \ge 0$ if $\successindices = \emptyset$).
All $t\ge t_0$ then must be model-improving, acceptable or inacceptable.
In all cases, the trust region radius $\radiust$ is never increased.
Due to \cref{assumption:surrogates}, the number of successive model-improvement 
steps is bounded above by $\maximprovements \in \mathbb{N}$.
Hence, $\radiust$ is decreased by a factor of
$\gamma\in [\gammasmallest, \gammasmall] \subseteq (0,1)$ at 
least once every $\maximprovements$ iterations.
Thus,
$$
	\sum_{t > t_0}^\infty \radiust
	\le
		N \sum_{i = 1}^\infty \gammasmall^i \radius\ita{t_0}
	=
		\frac{N\gammasmall}{1-\gammasmall}\radius\ita{t_0},
$$
and $\radiust$ \textbf{must go to zero} for $t\to \infty$.

Clearly, for any $\tau \ge t_0$, the iterates 
(and trust region centers) $\ve x\ita{\tau}$ and $\ve x\ita{t_0}$ 
cannot be further apart than the sum of all subsequent 
trust region radii, i.e.,
$$
	\norm{ \ve x\ita{\tau} -\ve x\ita{t_0} }
	\le
		\sum_{t\ge t_0}^\infty \radiust
	\le
		\frac{N\gammasmall}{1-\gammasmall}\radius\ita{t_0}.
$$
The RHS goes to zero as we let $t_0$ go to infinity and so 
must the norm on the LHS, i.e.,
\begin{equation}
 	\lim_{t_0 \to \infty} \norm{ \ve x\ita{\tau} -\ve x\ita{t_0} } = 0.
	\label{eqn:finite_success_distance_to_zero}
\end{equation}
Now let $\tau = \tau(t_0) \ge t_0$ be the first iteration 
index so that $\ve m\ita{\tau}$ is fully linear.
Then
$$
 \left| \modshortmit{t_0} \right|
 \le
 	\left|
		\modshortit{t_0}
			- \modshortit{\tau}
	\right|
	+
	\left|
		\modshortit{\tau}
			- \modshortmit{\tau}
	\right|
	+
	\left| \modshortmit{\tau} \right|
$$
and for the terms on the right and for $t_0\to \infty$, we find:
\begin{itemize}
	\item 
		Because of \cref{assumption:objective_gradients_lipschitz,assumption:feasible_set_compact_convex} 
		and \cref{theorem:omega_uniformly_continuous} $\modomega{\bullet}$ is Cauchy-continuous 
		and with \eqref{eqn:finite_success_distance_to_zero} the first term goes to zero.
	\item 
		Due to \cref{theorem:omegas_fully_linear} the second term 
		is in $\mathcal O(\radius\ita{\tau})$ and goes to zero.
	\item 
		Suppose the third term does not go to zero as well, i.e., 
		$\{\modomegamt{\ve x\ita{\tau}}\}$ is bounded below by a positive constant.
		Due to \cref{assumption:surrogates,assumption:feasible_set_compact_convex} the 
		iterates $x\ita{\tau}$ are not Pareto-critical for \eqref{eqn:mopmt} and
		because of $\radius\ita{\tau}\to 0$ and \cref{theorem:small_radius_very_successful} 
		there would be a successful iteration, a contradiction.
		Thus the third term must go to zero as well.
\end{itemize}
We conclude that the left side, $\modomega{\ve x\ita{t_0}}$, goes to zero as well for $t_0\to \infty$.
\end{proof}

We now address the case of infinitely many successful iterations,
first for the surrogate measure $\modomegamt{\bullet}$ and then for $\modomega{\bullet}$.
We show that the criticality measures are not bounded away from zero.\\
We start with the observation that in any case the trust region radius converges to zero:
\begin{Lemma} 
	If
	\cref{assumption:feasible_set_compact_convex,assumption:model_hessian_bounded_above,assumption:objective_gradients_lipschitz,assumption:surrogates,assumption:sufficient_decrease}
	hold, then the subsequence of trust region radii generated by \cref{algo:algorithm1} goes to zero, i.e.,
	$
		\lim_{t\to \infty} \radiust = 0.
	$
	\label{theorem:radius_to_zero}
\end{Lemma}

\begin{proof}
%See also \cite[Lemma 10.9]{conn:derivative_free_2009}.
We have shown in the proof of \cref{theorem:finitely_many_acceptable_iterations_convergence} that this is the case for
finitely many successful iterations.

Suppose there are infinitely many successful iterations.
Take any successful index $t\in \successindices$.
Then $\rho\itat \ge \nusuccess$ and
from \cref{assumption:sufficient_decrease}
%by the interpolation properties of $\ve m$ at all iterates
it follows for $\ve x\ita{t+1} = \xtrialt = \xitat + \stept$ that 
%\spfoot{Vorschlag: Eine Zeile und keine Nummber. Allgemein am besten alle Nummern vermeiden, die nicht referenziert werden.}
\begin{align}
	\Phif(\xitat) - \Phif(\xtrialt)
	&\ge
		\nusuccess \left(
			\Phimt(\xitat) - \Phimt(\xtrialt)
			\right)
	\stackrel{\mathclap{\eqref{eqn:sufficient_decrease_modified_prototype}}}
	\ge
		\nusuccess \sufficientconstantm \modshortmit{t} \min
		\left\{
			\frac{\modshortmit{t}}{\normconst\hessboundm},
			\radiust
		\right\}. \notag
\shortintertext{The criticality step ensures that $\modshortmit{t}\ge \min\left\{\epscrit, \dfrac{\radiust}{\mu}\right\}$ so that}
	\Phif(\xitat) - \Phif(\xtrialt)
	&\ge
		\nusuccess \sufficientconstantm \min\left\{\epscrit, \dfrac{\radiust}{\mu}\right\} \min
		\left\{
			\frac{\radiust}{\mu\normconst\hessboundm},
			\radiust
		\right\}\ge 0.
	\label{eqn:radius_to_zero}
\end{align}
Now the right hand side has go to zero: Suppose it was bounded below by a positive constant $\varepsilon>0$.
We could then compute a lower bound on the improvement from the first iteration with index $0$ up to $t+1$ by summation
$$
	\Phif(\ve x\ita 0) - \Phif (\ve x\ita{t+1})
	\ge
	\sum_{\tau\in \successindices_t}
	\Phif(\ve x\ita{\tau}) - \Phif(\ve x\ita{\tau+1})
	\ge
	\left|
		\successindices_t
	\right|
	\varepsilon
$$
where $\successindices_t = \successindices \cap\{0,\ldots,t\}$ are all successful 
indices with a maximum index of $t$.
Because $\successindices$ is unbounded, the right side diverges for $t\to \infty$ 
and so must the left side in contradiction to $\Phif$ being bounded
below by \cref{assumption:objectives_bounded}.
From \eqref{eqn:radius_to_zero} we see that this implies $\radiust \to 0$ for $t \in \successindices, t\to \infty$.\\
Now consider any sequence $\mathcal T \subseteq N$ of indices that are not necessarily successful, 
i.e., $\left|\mathcal T \setminus \successindices\right| \ge 0$.
The radius is only ever increased in successful iterations and at most by a factor of $\gammabig$. 
Since $\successindices$ is unbounded, there is for any $\tau\in \mathcal T$ a largest $t_\tau\in \successindices$ with $t_\tau \le \tau$.
Then $\radius\ita \tau \le \gammabig\radius\ita{t_\tau}$ and because of $\radius\ita{t_\tau} \to 0$ it follows that
$$
\lim_{\substack{\tau\in \mathcal T,\\ \tau\to \infty}} \radius\ita{\tau} = 0,
$$
which concludes the proof.
\end{proof}

\begin{Lemma}
	Suppose 		
	\cref{assumption:model_hessian_bounded_above,assumption:objective_gradients_lipschitz,assumption:surrogates,assumption:sufficient_decrease,assumption:feasible_set_compact_convex,assumption:objectives_bounded} 
	hold.
	For the iterates produced by \cref{algo:algorithm1} it holds that
	$$
		\liminf_{t\to \infty} \modomegamt{ \xitat } = 0.
	$$
	\label{theorem:convergence_liminf}
\end{Lemma}

\begin{proof}
%See also \cite[Lemma 10.10]{conn:derivative_free_2009},\cite[Lemmata 5.6]{conn:trm_framework} and \cite[Lemma 4.20]{thomann}.\\
For a contradiction, suppose that
$
	\liminf_{t\to \infty} \modomegamt{ \xitat } \ne 0.
$
Then there is a constant $\lbomega > 0$ with $\modshortmit{t}\ge \lbomega$ for all $t\in \mathbb{N}_0$.
According to \cref{lemma:lbomega_lbradius}, there exists a constant $\lbradius>0$ with $\radiust\ge \lbradius$ for all $t$.
This contradicts \cref{theorem:radius_to_zero}.
\end{proof}

The next result allows us to transfer the result to $\modomega{\bullet}$.
\begin{Lemma}
	Suppose 
	\cref{assumption:feasible_set_compact_convex,assumption:objective_gradients_lipschitz,assumption:surrogates} 
	hold.
	For any subsequence $\{t_i\}_{i\in \mathbb{N}}\subseteq \mathbb{N}_0$ of
	iteration indices of \cref{algo:algorithm1} with
	\begin{equation}
		\lim_{i\to \infty} \modomegamit{t_i}
			%\ve x \ita{t_i} }
		 = 0,
		\label{eqn:subsequence_to_zero}
	\end{equation}
	it also holds that
	\begin{equation}
		\lim_{i\to \infty} \modomegait{t_i} = 0. %{\ve x\ita{t_i}} = 0.
		\label{eqn:subsequence_to_zero2}
	\end{equation}
	\label{theorem:subsequence_to_zero}
\end{Lemma}

\begin{proof}
	%See also \cite[Lemma 5.7]{conn:trm_framework}.

	By \eqref{eqn:subsequence_to_zero}, $\modshortmit{t_i} < \epscrit$ for sufficiently large $i$.
	If $\ve x \ita{t_i}$ is critical for \eqref{eqn:mop}, then the result follows from \cref{theorem:inifinite_criticality_loop}.
	Otherwise, $\ve m\ita{t_i}$ is fully linear on $\ball{\ve x \ita{t_i}}{\radius\ita{t_i}}$ for some $\radius\ita{t_i}\le \mu\modshortmit{t_i}$.
	From \cref{theorem:omegas_fully_linear} it follows that
	$$
	\left|
		\modshortmit{t_i} - \modshortit{t_i}
	\right|
	\le \sqrt{c} \fullylinearconstantdf \radius\ita{t_i}
	\le \sqrt{c} \fullylinearconstantdf \mu \modshortmit{t_i}.
	$$
	The triangle inequality yields
	$$
	 \modshortit{t_i}
		\le
			\left|
				 \modshortit{t_i}
				- \modshortmit{t_i}
			\right|
			+
			\modshortmit{t_i}
		\le
		(\sqrt{c} \fullylinearconstantdf \mu + 1 )\modshortmit{t_i}
	$$
	for sufficiently large $i$ and \eqref{eqn:subsequence_to_zero} then implies \eqref{eqn:subsequence_to_zero2}.
\end{proof}

The next global convergence result immediately follows from \cref{theorem:finitely_many_acceptable_iterations_convergence,theorem:convergence_liminf,theorem:subsequence_to_zero}:
\begin{Theorem}
	Suppose 		
	\cref{assumption:feasible_set_compact_convex,assumption:model_hessian_bounded_above,assumption:sufficient_decrease,assumption:objective_gradients_lipschitz,assumption:surrogates,assumption:sufficient_decrease,assumption:objectives_bounded} 
	hold.
	Then
	$
		\liminf_{t\to \infty}
		\modomega{\xitat} = 0.
	$
	\label{theorem:convergence_liminf_true}
\end{Theorem}

This shows that if the iterates are bounded, then there is a subsequence of iterates in $\mathbb{R}^n$ approximating a
Pareto-critical point.
We next show that \emph{all} limit points of a sequence generated by \cref{algo:algorithm1} are Pareto-critical. 
%\spfoot{Ich fänd es super, wenn wir das schon am Ende von Kap.\ 2 erwähnen könnten, also direkt nach dem Alg.\ 1. In etwa so etwas wie ``The main theoretical result of this paper will be to prove convergence for a large class of surrogate models: \textbf{Theorem 1:} ... ''. Dann kann ein Hinweis folgen, dass die Konvergenzanalyse und der Beweis in Kap.\ 5 kommen. Den Verweis könnte man auch nach der Einführung von Alg.\ 2 nochmal machen. Ggf.\ könnte man auch ein Theorem 1 einführen, das sehr viel weniger formal ist (so etwas wie ``under some assumptions that are common in optimization as well as surrogate modeling, Algorithms 1 and 2 converge to Pareto critical points.'') und dann sagen, dass in Kap.\ 5 die formale Variante definiert und auch bewiesen wird.}

\begin{Theorem}
	Suppose
	\cref{assumption:feasible_set_compact_convex,assumption:model_hessian_bounded_above,assumption:sufficient_decrease,assumption:model_hessian_bounded_above,assumption:omegaconstant,assumption:objective_gradients_lipschitz,assumption:surrogates,assumption:objectives_bounded} 
	hold.
	Then
	$
		\lim_{t\to \infty}
		\modomega{\xitat} = 0.
	$
\end{Theorem}
\begin{proof}
%	See \cite[Lemma 10.13]{conn:derivative_free_2009},\cite[Lemma 5.9]{conn:trm_framework} and \cite[Th. 4.21]{thomann}.\\
	We have already proven the result for finitely many successful iterations, see \cref{theorem:finitely_many_acceptable_iterations_convergence}.
	We thus suppose that $\successindices$ is unbounded.

	For the purpose of establishing a contradiction, suppose that there exists a sequence
	$\left\{ t_j\right\}_{j\in \mathbb{N}}$ of indices that are
	successful or acceptable with
	\begin{equation}
		\modshortit{t_j} \ge 2\varepsilon > 0
		\quad \text{for some $\varepsilon>0$ and all $j$.}
		\label{eqn:two_epsilon_t_j}
	\end{equation}
	We can ignore model-improving and inacceptable iterations:
	During those the iterate does not change and we find a larger acceptable or successful index 
	with the same criticality value.

	From \cref{theorem:convergence_liminf_true} we obtain that for every such $t_j$, there exists
	a first index $\tau_j>t_j$ such that $\modomega{\ve x\ita{\tau_j}} < \varepsilon$.
	We thus find another subsequence indexed by $\{\tau_j\}$ such that
	\begin{equation}
		\modshortit{t} \ge \varepsilon
		\text{ for $t_j\le t<\tau_j$ and }
		\modshortit{\tau_j} < \varepsilon.
		\label{eqn:subsequence_definitions}
	\end{equation}	
	Using \eqref{eqn:two_epsilon_t_j} and \eqref{eqn:subsequence_definitions}, it also follows from a triangle inequality that
	\begin{equation}
		\left|
			\modshortit{t_j} - \modshortit{\tau_j}
		\right|
		\ge
			\modshortit{t_j}
			-
			\modshortit{\tau_j}
		>
			2\varepsilon
			-
			\varepsilon = \varepsilon  \qquad \forall j\in \mathbb{N}.
		\label{eqn:distance_bigger_eps}
	\end{equation}
	With $\{t_j\}$ and $\{\tau_j\}$ as in \eqref{eqn:subsequence_definitions}, define the following subset
	set of indices
	$$
		\mathcal T =
		\left\{
			t \in \mathbb{N}_0: \exists j\in \mathbb{N} \text{ such that } t_j \le t < \tau_j
		\right\}.
	$$
	By \eqref{eqn:subsequence_definitions} we have $\modshortit{t}\ge \varepsilon$ for $t\in \mathcal T$, and
	due to \cref{theorem:subsequence_to_zero}, we also know that then $\modshortmit{t}$ cannot
	go to zero neither, i.e., there is some $\varepsilon\subm>0$ such that
	\begin{equation*}
	 	\modshortmit{t} \ge
		\varepsilon\subm > 0 \qquad \forall t \in \mathcal T.
	\end{equation*}
	From \cref{theorem:radius_to_zero} we know that $\radiust \xrightarrow{t\to \infty} 0$ so that by
	\cref{theorem:small_radius_very_successful_true_prolem}, any sufficiently large $t\in \mathcal T$ must be either
	successful or model-improving (if $\vemt$ is not fully linear).
	For $t\in \mathcal T\cap \successindices$, it follows from \cref{assumption:sufficient_decrease} that
	\begin{equation*}
		\Phif( \ve x\ita t ) - \Phif( \ve x\ita {t+1} )
		\ge
		\nusuccess \left(
			\Phim( \ve x\ita t ) - \Phim( \ve x\ita {t+1} )
		\right)
		\ge
		\nusuccess \sufficientconstantm \varepsilon\subm \min
		\left\{
			\frac{\varepsilon\subm}{\normconst \hessboundm},
			\radiust
		\right\}
		\ge 0.
		\label{eqn:sufficient_decrease_main_theorem}
	\end{equation*}
	%Hence $\left\{\Phif(\ve x\ita t)_{t\in{\tilde{\acceptableindices}}}\right\}$ is 
	%monotonically decreasing (for large $t$) and by \cref{assumption:objectives_bounded} bounded below.
	%The sequence is thus convergent and
	%$$\lim_{\substack{t\to \infty\\t\in{\tilde{\acceptableindices}}}}\Phif( \ve x\ita t ) - \Phif( \ve x\ita {t+1} ) =0$$
	%which in turn implies $\radiust \to 0$.
	If $t\in \mathcal T\cap \successindices$ is sufficiently large, we have
	$\radiust \le \dfrac{\varepsilon\subm}{\normconst \hessboundm}$ and
	$$
		\radiust \le \frac{1}{\nusuccess \sufficientconstantm \varepsilon\subm}
		\left(
			\Phif( \ve x\ita t ) - \Phif( \ve x\ita {t+1} )
		\right).
	$$
	Since the iteration is either successful or
	model-improving for sufficiently large $t\in \mathcal T$, 
	and since $\xitat = \ve x\ita{t+1}$ for a model-improving iteration,
	we deduce from the previous inequality that
	$$
	\norm{ \ve x\ita{t_j} - \ve x\ita{\tau_j} }
	\le
		\sum_{\substack{t = t_j,\\t\in \mathcal T \cap \successindices}}^{\tau_j -1}
		\norm{ \ve x\ita t - \ve x\ita{t +1} }
	\le
		\sum_{\substack{t = t_j,\\t\in \mathcal T \cap \successindices}}^{\tau_j -1} \radiust
	\le
		\frac{1}{\nusuccess\sufficientconstantm\varepsilon\subm}
		\left(
			\Phif( \ve x\ita {t_j} ) - \Phif( \ve x\ita {\tau_j} )
		\right)
	$$
	for $j\in \mathbb{N}$ sufficiently large.
	The sequence $\left\{ \Phif(\xitat)\right\}_{t\in \mathbb{N}_0}$ is bounded below (\cref{assumption:objectives_bounded}) and
	monotonically decreasing by construction.
	Hence, the RHS above must converge to zero for $j\to \infty$.
	This implies $\lim_{j\to \infty}\norm{ \ve x\ita{t_j} - \ve x\ita{\tau_j} } = 0$.\\
	Because of \cref{assumption:feasible_set_compact_convex,assumption:objective_gradients_lipschitz}, $\modomega{\bullet}$
	is uniformly continuous so that then
	$$
		\lim_{j\to \infty}
			\modomega{\ve x\ita{t_j}}
			-\modomega{\ve x\ita{\tau_j}}
		=0,
	$$
	which is a contradiction to \eqref{eqn:distance_bigger_eps}.
	Thus, no subsequence of acceptable or successful indices as in \eqref{eqn:two_epsilon_t_j} can exist.

	% There must not be such a sequence even if we allow inacceptable indices.
	% %i.e., $\{t_j\}_{j\in \mathbb N} \setminus \acceptableindices \ne \emptyset$.
	% This is again because $\successindices$ (and thereby $\acceptableindices$) 
	% is unbounded and the sequence
	% $$
	% 	\left\{ 
	% 		z_j \in \mathbb N_0 : z_j = 
	% 		\min_{z\in \acceptableindices} z \text{ s.t. $z \ge t_j$ }
	% 	\right\} \subseteq \acceptableindices
	% $$
	% would satisfy $\modshortit{z_j} = \modshortit{t_j}$ for all $j$ (because $\xitat$ 
	% does not change in inacceptable iterations) 
	% and	result in the same contradiction as above.

	% All in all, we have shown that no subsequence of iterates can exist that is bounded below by a positive constant.
	% From \cref{theorem:convergence_liminf_true} it follows that $\lim_{t\to \infty} \modshortit{t} = 0.$
\end{proof}

\section{Numerical Examples}
\label{section:examples}

In this section we provide some more details on the 
actual implementation of \cref{algo:algorithm1} and
present the results of various experiments.
We compare different surrogate model types with regard
to their efficacy (in terms of expensive objective evaluations) 
and their ability to find Pareto-critical points.

\subsection{Implementation Details}
We implemented the algorithm in the Julia language.
The \verb|OSQP| solver \cite{osqp} was used to solve \eqref{eqn:descent_direction_problemm}.
For non-linear problems we used the \verb|NLopt.jl| \cite{nlopt} package. 
More specifically we used the \verb|BOBYQA| algorithm \cite{bobyqa} in conjunction with 
\verb|DynamicPolynomials.jl| \cite{DynamicPolynomials} for the Lagrange polynomials 
and the population based \verb|ISRES| method \cite{isres} for the 
Pascoletti-Serafini subproblems.
The derivatives of cheap objective functions were obtained by means 
of automatic differentiation \cite{forwarddiff} and Taylor models used 
\verb|FiniteDiff.jl|.%\cite{finitediff}.

In accordance with \cref{algo:algorithm1} we perform the 
shrinking trust region update via
$$
\radiusnext \leftarrow 
\begin{cases}
    \gammasmallest \radiust&\text{if $\rho\itat < \nuaccept$,}\\
    \gammasmall \radiust&\text{if $\rho \itat < \nusuccess$.}
\end{cases}
$$
Note that for box-constrained problems we internally scale the feasible set
to the unit hypercube $[0,1]^n$ and all radii are measured with regard
to this scaled domain.

For \textbf{stopping} we use a combination of different criteria:
\begin{itemize}
    \item We have an upper bound $\maxiter \in \mathbb N$ on the maximum number of iterations 
    and an upper bound $\maxexpensive \in \mathbb N$ on the number of expensive objective evaluations.
    \item The surrogate criticality naturally allows for a stopping test and due to \cref{lemma:lbomega_lbradius}
    the trust region radius can also be used (see also \cite[Sec. 5]{thomann}).
    We combine this with a relative tolerance test and stop if 
    $$
        \radiust \le \radiusmin \text{ OR }
        \left( \radiust \le \radiuscrit \text{ AND } \norm{\ve s\itat} \le \stepsmin\right).
    $$
    \item At a truly critical point the criticality loop \criticalityRoutine runs infinitely. 
    We stop after a maximum number $\maxcritloops \in \mathbb N_0$ of iterations.
    If $\maxcritloops$ equals 0 the algorithm effectively stops for small $\modshortmit{t}$ values.
\end{itemize}

\subsection{A First Example}

We tested our method on a multitude of academic test
problems with a varying number of decision variables 
$n$ and objective functions $k$.
We were able to approximate Pareto-critical points in both cases,
if we treat the problems as heterogenous and if we declare them as
expensive.
We benchmarked RBF against polynomial models, because in \cite{thomann}
it was shown that a trust region method using second degree Lagrange 
polynomials outperforms commercial solvers on scalarized problems.
Most often, RBF surrogates outperform other model types with regard to the 
number of expensive function evaluations.

\begin{figure}[H]
    \centering 
    \includegraphics[width=\linewidth]{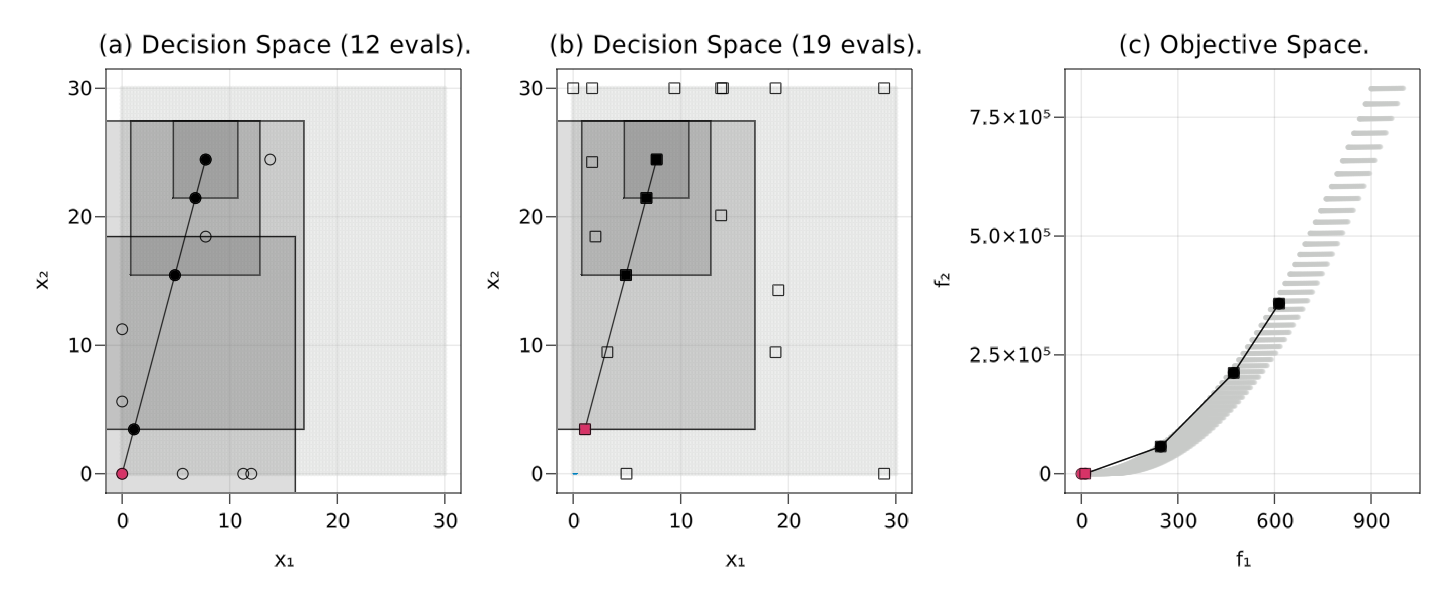}
    \caption{
        Two runs with maximum number of expensive evaluations set to 20 (soft limit). %
        Test points are light-gray, the iterates are black, final iterate is red, white markers show other points where the objectives are evaluated. %
        The successive trust regions are also shown. %
        \textbf{(a)} Using RBF surrogate models we converge to the optimum using only 12 expensive evaluations. %
        \textbf{(b)} Quadratic Lagrange models do not reach the optimum using 19 evaluations. %
        \textbf{(c)} Iterations and test points in the objective space. %
    }
    \label{fig:mht}
\end{figure}

This is illustrated in \cref{fig:mht}. 
It shows two runs of \cref{algo:algorithm1} on the non-convex 
problem \eqref{eqn:mht}, taken from \cite{thomann_diss}:
\begin{equation}
    \begin{aligned}
        \min_{\ve x\in \feas}
        \begin{bmatrix}
            x_1 + \ln (x_1) + x_2^2,\\
            x_1^2 + x_2^4
        \end{bmatrix}, \, \feas = [\varepsilon,30]\times [0,30] \subseteq \mathbb R^2, \varepsilon = 10^{-12}.
    \end{aligned}
    \tag{T6}
    \label{eqn:mht}
\end{equation}

The first objective function is treated as expensive while the second is cheap.
The only Pareto-optimal point of \eqref{eqn:mht} is $[\varepsilon, 0]^T$.
When we set a very restrictive limit of $\maxexpensive = 20$ then we run 
out of budget with second degree Lagrange surrogates before we reach the optimum, 
see \cref{fig:mht} \textbf{(b)}.
As evident in \cref{fig:mht} \textbf{(a)}, surrogates based on 
(cubic) RBF do require significantly less training data. 
For the RBF models the algorithm stopped after $2$ critical loops
and the model refinement during these loops is made clear by the samples
on the problem boundary converging to zero.
The complete set of relevant parameters for the test runs is given in \cref{table:params1}.
We used a strict acceptance test and the strict Pareto-Cauchy step.

\begin{specialtable}[htb]
    \caption{Parameters for \cref{fig:mht}, radii relative to $[0,1]^n$.}
    \label{table:params1}
\begingroup
\small 
\setlength{\tabcolsep}{0.5em}
\begin{tabular}{ccccccccccccc}
    \toprule 
    $\epscrit$ & $\maxexpensive$ & $\maxcritloops$ & $\mu$ & $\beta$ 
    & $\radiusmax$ & $\radiusmin$ & $\radius\ita{0}$ & $\nuaccept$ & $\nusuccess$     
    & $\gammasmallest$ & $\gammasmall$ & $\gammabig$
    \\
    \midrule 
    $\ten{-3}$ & $20$ & $2$ & $2\cdot\ten{3}$ & $\ten{3}$
    & $0.5$ & $\ten{-3}$ & $0.1$ & $0.1$ & $0.4$ 
    & $0.51$ & $0.75$ & $2$
    \\
    \bottomrule
\end{tabular}
\endgroup
\end{specialtable}

\subsection{Benchmarks on Scalable Test-Problems}

\begin{figure}[H]
    \centering
    \includegraphics[width=\linewidth]{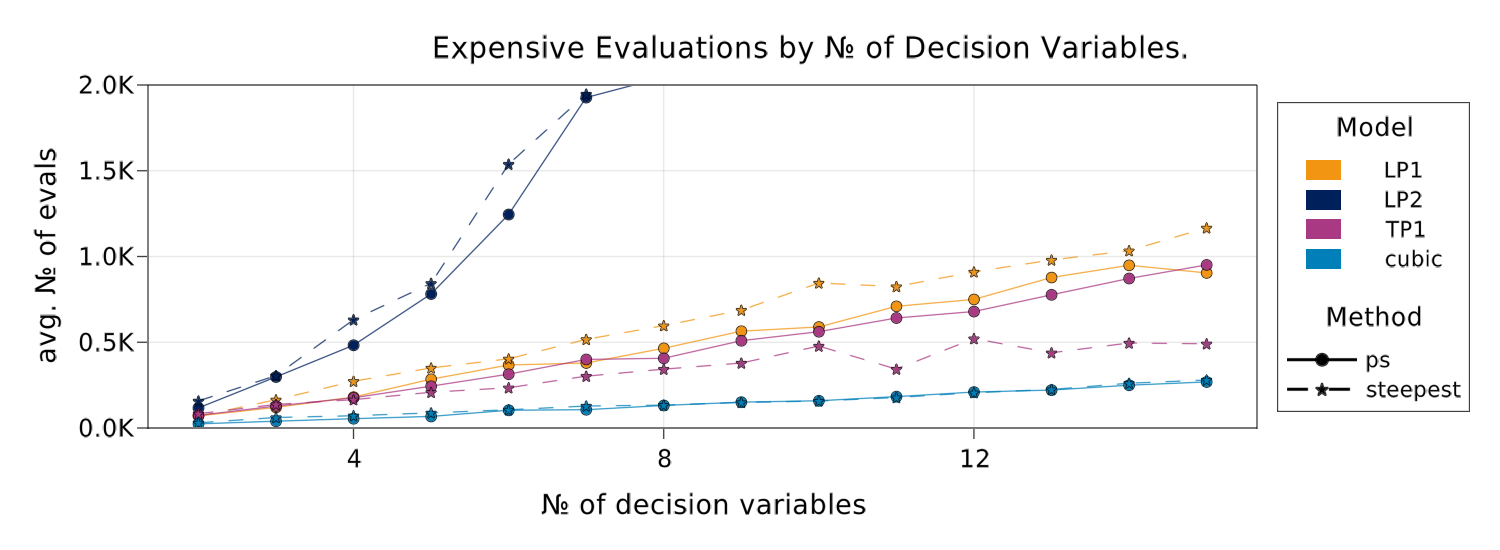}
    \caption{Average number of expensive objective evaluations by number of decision variables $n$, surrogate type and descent method. %
    LP1 are Linear Lagrange models, LP2 quadratic Lagrange models, TP1 are linear Taylor polynomials based on 
    finite differences and \emph{cubic} refers to cubic RBF models. %
    Steepest descent and Pascoletti-Serafini were tested on scalable problems, and 12 runs were performed per setting.}
    \label{fig:avg_evals_per_dim}
\end{figure}

To assess the performance with a growing number of decision variables $n$, we 
performed tests on scalable problems of the ZDT and DTLZ family \cite{zdt,dtlz}.
\cref{fig:avg_evals_per_dim} shows results for the bi-objective problems ZDT1-ZDT3 
and for the $k$-objective problems DTLZ1 and DTLZ6 
(we used $k = \max\{2, n-4\}$ objectives).
All problems are box constrained. 
Eight feasible starting points were generated for each problem setting, 
i.e., for each combination of $n$, a test problem 
and a descent method).

In all cases the first objective was considered cheap and all other objectives expensive.
First and second degree Lagrange models were compared against linear Taylor models and 
(cubic) RBF surrogates.
The Lagrange models were built using a $\Lambda$-poised set, with $\Lambda = 1.5$.
In the case of quadratic models we used a precomputed set of points for $n\ge 6$.
The Taylor models used finite differences and points outside of box constraints were
simply projected back onto the boundary.
The RBF models were allowed to include up to $(n+1)(n+2)/2$ training 
points from the database if $n\le 10$ and else the maximum number of points was $2n+1$.
All other parameters are listed in \cref{table:params2}.

\begin{specialtable}[H]
    \caption{Parameters for \cref{fig:avg_evals_per_dim}, radii relative to $[0,1]^n$.
    $\theta_1$ is used to construct Lagrange and RBF models in an enlarged trust region, 
    $\theta_2$ is used only for RBF, see \cref{section:rbf}.}
    \label{table:params2}
\begingroup
\small 
\setlength{\tabcolsep}{0.5em}
\begin{tabular}{c|ccccccccccccccc}
    \toprule 
    \textbf{Parameter} &
    $\epscrit$ & $\maxiter$ & $\maxexpensive$ & $\maxcritloops$ & $\mu$ & $\beta$ 
    & $\radiusmax$ & $\radiuscrit$ & $\stepsmin$ & $\radiusmin$ & $\radius\ita{0}$ 
    %& $\nuaccept$ & $\nusuccess$     
    %& $\gammasmallest$ & $\gammasmall$ & $\gammabig$
    \\
    \midrule 
    \textbf{Value} &
    $\ten{-3}$ & $100$ & $n^2\cdot \ten{3}$ & $10$ & $2\cdot\ten{3}$ & $\ten{3}$
    & $0.5$ & $\ten{-3}$ & $\ten{-8}$ & $\ten{-6}$ & $0.1$ 
    %& $0.0$ & $0.4$ 
    %& $0.51$ & $0.75$ & $2$
    \\
    \bottomrule
    \toprule 
    \textbf{Parameter} &
    $\nuaccept$ & $\nusuccess$     
    & $\gammasmallest$ & $\gammasmall$ & $\gammabig$ & $\theta_2$ & $\theta_1$
    \\
    \midrule 
    \textbf{Value} &
    $0.0$ & $0.4$ 
    & $0.51$ & $0.75$ & $2$ & $5$ & $2$ 
    \\
    \bottomrule
\end{tabular}
\endgroup
\end{specialtable}

As expected, the second degree Lagrange polynomials require the most 
objective evaluations and the quadratic dependence on $n$
is clearly visible in \cref{fig:avg_evals_per_dim}, and
the quadratic growth of the dark-blue line continues for $n\ge 8$. %, the y-Axis is cropped
%so as not to squash the other curves.
On average, the Taylor models perform better than the linear Lagrange polynomials -- despite
requiring more evaluations per iteration.
This is possibly due to more accurate derivative information and resulting faster convergence.
The Lagrange models do slightly better when the Pascoletti-Serafini step calculation
is used (see \cref{appendix:section_steps}).
By far the least evaluations are needed for the RBF models: The light-blue line
consistently stays below all other data points. 
Often, the average number of evaluations is less than half that of the second best method.

\begin{figure}[H]
    \centering
    \includegraphics[width=.49\linewidth]{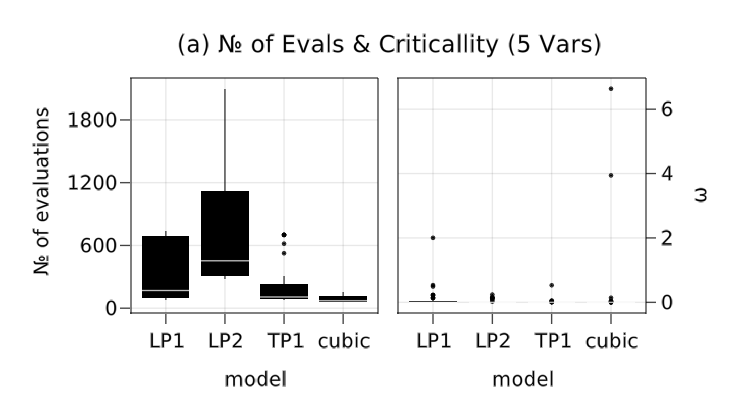}
    \includegraphics[width=.49\linewidth]{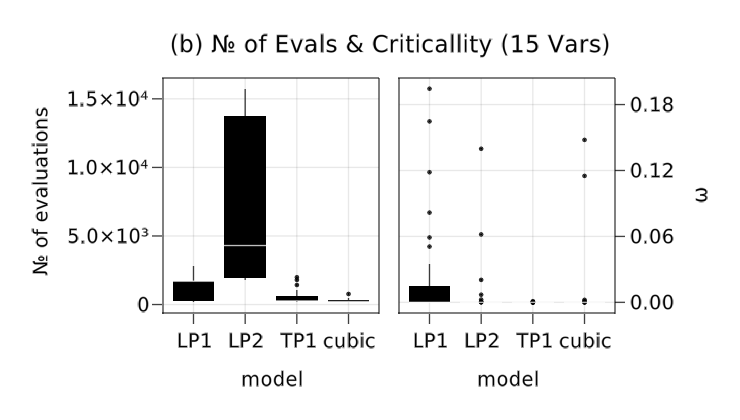}
    \caption{Box-plots of the number of evaluations and the solution criticality for $n=5$ and 
    $n=15$ for the steepest-descent runs from \cref{fig:avg_evals_per_dim}.}
    \label{fig:boxplots}
\end{figure}

\begin{figure}[H]
    \centering
    \includegraphics[width = .99\linewidth]{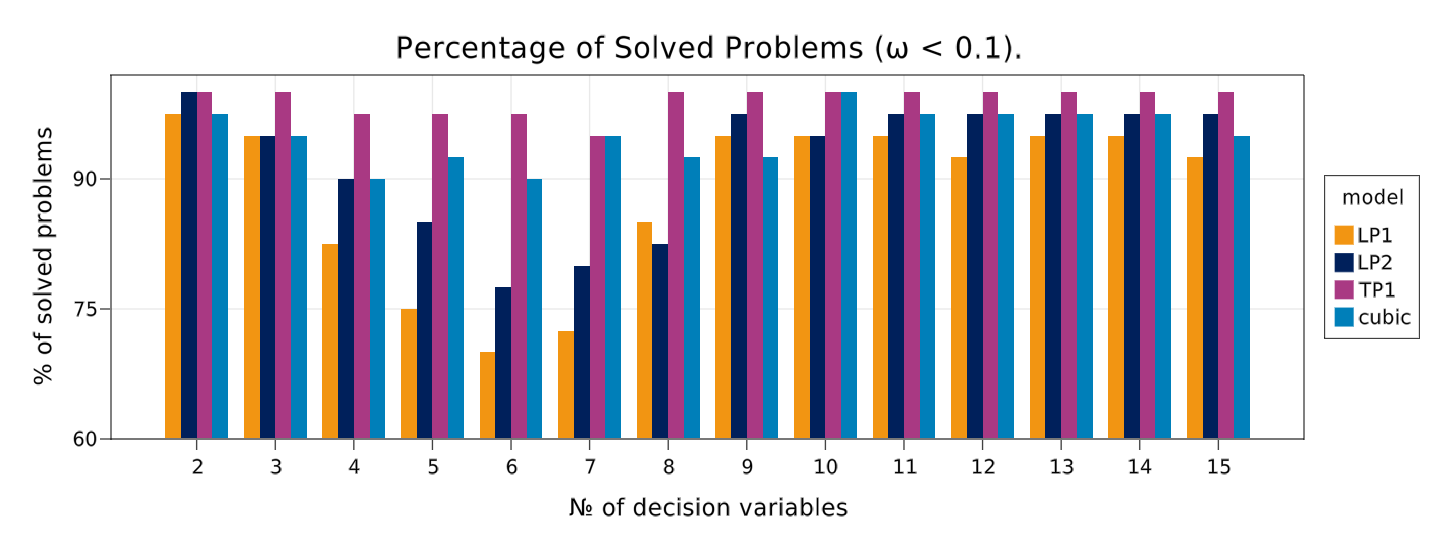} 
    \caption{Percentage of solved problem instances, i.e., test runs were the final solution 
    criticality has a value below 0.1.
    Per model and $n$-value there were 40 runs.}   
    \label{fig:omega_nvars_models}
\end{figure}

\cref{fig:boxplots} illustrates that not only do RBF perform better on average, 
but also overall.
With regards to the final solution criticality, there are a few outliers when 
the method did not converge using RBF models.
However, in most cases the solution criticality is acceptable, see \cref{fig:boxplots} \textbf{(b)}.
Moreover, \cref{fig:omega_nvars_models} shows that a good percentage of problem instances 
is solved with RBF, especially when compared to linear Lagrange polynomial models.
Note, that in cases where the true objectives are not differentiable at the final iterate,
$\omega$ was set to 0 because the selected problems are non-differentiable only in 
Pareto-optimal points.

\begin{figure}[H]
    \centering
    \includegraphics[width = .99\linewidth]{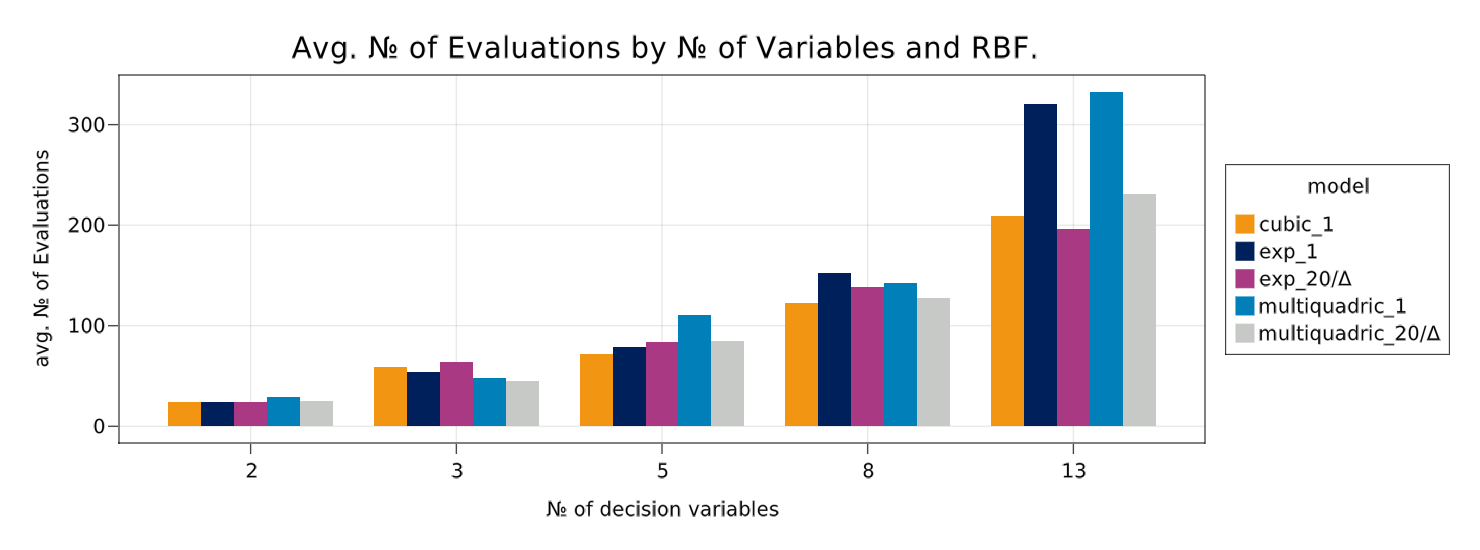}    
    \caption{Influence of a adaptive shape radius on the performance of
    RBF models (tested on ZDT3).}
    \label{fig:rbf_shape_params}
\end{figure}

Furthermore, we compared the RBF kernels from \cref{table:rbf}.
In \cite{wild:orbit}, the cubic kernel performs best on single-objective problems while the Gaussian
does worst.
As can be seen in \cref{fig:rbf_shape_params} this holds for multiple objective 
functions, too. 
The dark-blue and the light-blue bars show that both the Gaussian and the Multiquadric 
require more function evaluations, especially in higher dimensions.
If, however, we use a very simple adaptive strategy to fine-tune the shape parameter, 
then both kernels can finish significantly faster.
The pink and the gray bar illustrate this fact. 
In both cases, the shape parameter was set to $\alpha = 20/\radiust$ in each iteration.
Nevertheless, cubic function (orange) appears to be a good choice in general.

\section{Conclusion}
\label{section:conclusion}

We have developed a trust region framework for
heterogeneous and expensive multiobjective optimization problems.
It is based on similar work \cite{qu,villacorta,ryu_derivative-free_2014,thomann}
and our main contributions are the 
integration of constraints and of radial basis function surrogates.
Subsequently, our method is is provably convergent for unconstrained problems and 
when the feasible set is convex and compact, while requiring significantly less 
expensive function evaluations due to a linear scaling of complexity
with respect to the number of decision variables.

For future work, several modifications and extensions can likely 
be transferred from the single-objective to the multiobjective case.
For examples, the trust region update can be made step-size-dependent (rather than $\rho\itat$ alone) to allow for 
a more precise model refinement, see \cite[Ch.\ 10]{conn_trust_region_2000}.
We have also experimented with the nonlinear CG method
\cite{nonlinear_cg} for a multiobjective Steihaug—Toint step 
\cite[Ch. 7]{conn_trust_region_2000} and early results look promising.

Going forward, we would like to apply our algorithm to a real world 
application, similar to what has been done in \cite{prinz_expensive_2020}.
Moreover, it would be desirable to obtain not just one but multiple
Pareto-critical solutions.
Because the Pascoletti-Serafini scalarization is still compatible with 
constraints, the iterations can be guided in image space by providing 
different global utopia vectors.
Furthermore, it is straightforward to use RBF with the heuristic 
methods from \cite{mht_front} for heterogeneous problems.
We believe that it should also be possible to propagate multiple 
solutions and combine the TRM method 
with non-dominance testing as has been done in the bi-objective 
case \cite{ryu_derivative-free_2014}.
One can think of other globalization strategies as well:
RBF models have been used in multiobjective Stochastic Search 
algorithms \cite{regis_multi-objective_2016} and trust region 
ideas have been included into population based strategies 
\cite{deb_trust_region_2019}.
It will thus bee interesting to see whether the theoretical convergence properties can be maintained
within these contexts by employing a careful trust-region management.
Finally, re-using the data sampled near the final iterate 
within a continuation framework like in \cite{schutze_pareto_2020} is a promising next step.

%%%%%%%%%%%%%%%%%%%%%%%%%%%%%%%%%%%%%%%%%%
\vspace{6pt} 

%%%%%%%%%%%%%%%%%%%%%%%%%%%%%%%%%%%%%%%%%%
%% optional
%\supplementary{The following are available online at \linksupplementary{s1}, Figure S1: title, Table S1: title, Video S1: title.}


%%%%%%%%%%%%%%%%%%%%%%%%%%%%%%%%%%%%%%%%%%
\authorcontributions{
    Conceptualization, M.B. and S.P.; methodology, M.B.; software, M.B.; validation, M.B. and S.P.; formal analysis, M.B. and S.P.; investigation, M.B.;
    writing---original draft preparation, M.B.; writing---review and editing, S.P.; visualization, M.B.; supervision, S.P.; 
    %project administration, S.P.; funding acquisition, S.P. 
    All authors have read and agreed to the published version of the manuscript.
}

\funding{
    This research has been funded by the European Union and the 
    German Federal State of North Rhine-Westphalia within the EFRE.NRW project ``SET CPS''.\\
    \includegraphics[width=.45\linewidth]{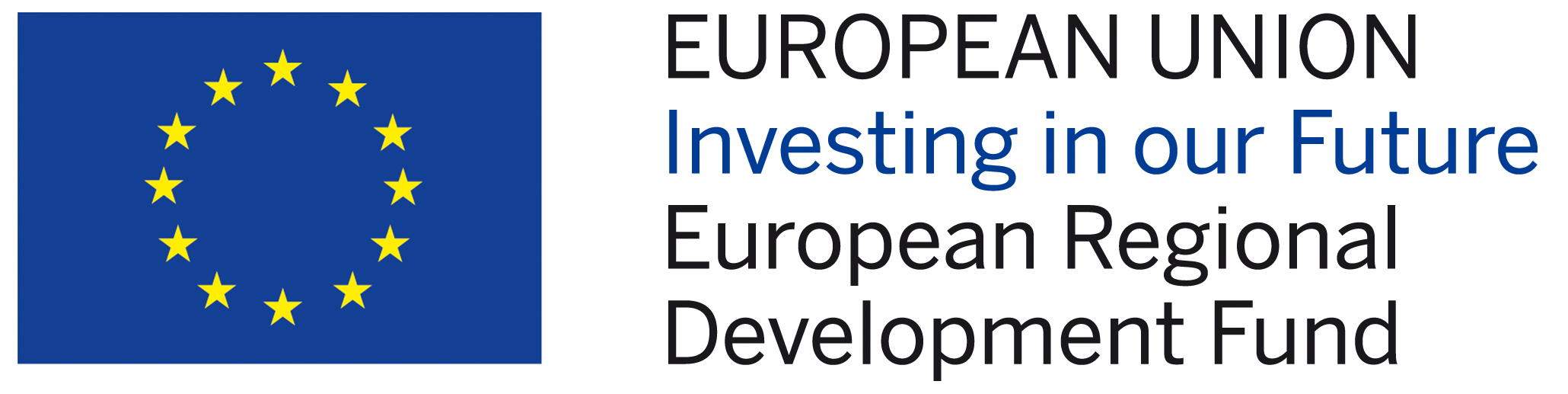}  
    \includegraphics[width =.45\linewidth]{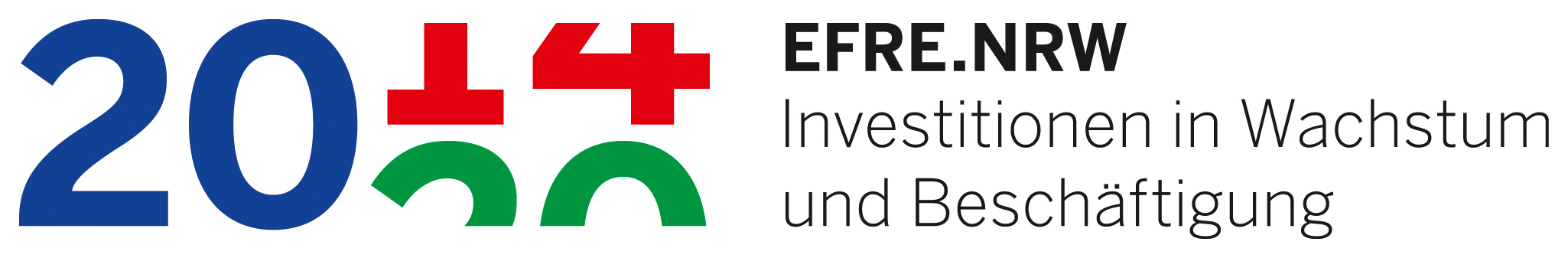}
}

% {\color{magenta}
% \acknowledgments{In this section you can acknowledge any support given which is not covered by the author contribution or funding sections. This may include administrative and technical support, or donations in kind (e.g., materials used for experiments).}
% }

\conflictsofinterest{``The authors declare no conflict of interest.''} 

%% Optional
% {\color{magenta}
% \sampleavailability{Samples of the compounds ... are available from the authors.}
% }

%%%%%%%%%%%%%%%%%%%%%%%%%%%%%%%%%%%%%%%%%%
%% Only for journal Encyclopedia
%\entrylink{The Link to this entry published on the encyclopedia platform.}

%%%%%%%%%%%%%%%%%%%%%%%%%%%%%%%%%%%%%%%%%%
%% Optional
% {\color{orange}
% \abbreviations{The following abbreviations are used in this manuscript:\\

% \noindent 
% \begin{tabular}{@{}ll}
%     TRM & Trust Region Method(s) \\
%     MO & Multiobjective Optimization \\
%     RBF & Radial Basis Function(s)
% \end{tabular}}
% }
% %%%%%%%%%%%%%%%%%%%%%%%%%%%%%%%%%%%%%%%%%%
%% Optional
\appendixtitles{yes} % Leave argument "no" if all appendix headings stay EMPTY (then no dot is printed after "Appendix A"). If the appendix sections contain a heading then change the argument to "yes".
\appendixstart
\appendix

\section{Miscellaneous Proofs}

\subsection{Continuity of the Constrained Optimal Value}
\label[Appendix]{appendix:omega_continuity}

In this subsection we show the continuity of $\omga{\ve x}$ 
in the constrained case, where $\omga{\ve x}$ is the negative
optimal value of \eqref{eqn:descent_direction_problem1}, i.e.,
\begin{align*}
    \omga{\ve x}:=
        &- \min_{ \ve d \in \feas - \ve x }  \max_{\ell=1,\ldots,k}
            \langle \gradf_\ell(\ve x), \ve d \rangle,
        \\
        &\text{s.t.}\;
        \norm{\ve d} \le 1.
\end{align*}
The proof of the continuity of $\omga{\ve x}$, as stated in \cref{theorem:p1_properties}, 
follows the reasoning from \cite{fukuda_drummond}, where continuity is shown for a 
related constrained descent direction program.

\begin{proof}[Proof of \cref{p1_properties2} in \cref{theorem:p1_properties}]
Let the requirements of \cref{p1_properties1} be fulfilled, i.e., 
let $\ve f$ be continuously differentiable and let 
$\feas \subset \mathbb{R}^n$ be convex and compact.
Further, let $\ve x$ be a point in $\feas$ and denote the minimizing direction in \eqref{eqn:descent_direction_problem1} by 
$\ve d(\ve x)$ and the optimal value by $\theta(\ve x)$. 
We show that $\theta(\ve x)$ is continuous, by which $\omga{\ve x} = - \theta(\ve x)$ is continuous as well.

\noindent First, note the following properties of the maximum function:
\begin{enumerate}
    \item \label[max_prop]{appendix:max_prop1} $\ve u \mapsto \max_\ell u_\ell$ is positively 
    homogenous and hence 
        $$ 
            \max_\ell \left( 
                \langle \gradf_\ell(\ve x), \ve d_1 \rangle 
                 + 
                \langle \gradf_\ell(\ve x), \ve d_2 \rangle 
                \right) 
            \le \max_\ell \langle \gradf_\ell(\ve x), \ve d_1 \rangle 
            + \max_\ell \langle \gradf_\ell(\ve x), \ve d_2 \rangle. 
        $$
    \item \label[max_prop]{appendix:max_prop2} 
        $\ve u \mapsto \max_\ell u_\ell$ is Lipschitz with constant 1 so that
        $$
        \abs{
            \max_\ell \langle \gradf(\ve x_1), \ve d_1 \rangle 
            - 
            \max_\ell \langle \gradf(\ve x_2), \ve d_2 \rangle
        }
        \le 
        \norm{
            \jacobianf_\ell(\ve x_1) \ve d_1 - \jacobianf_\ell(\ve x_2) \ve d_2
        },
        $$
        for both the maximum and the Euclidean norm.
\end{enumerate}
Now let $\{ \ve x\ita{t} \}\subseteq \feas$ be a sequence with $\ve x\ita{t}\to \ve x$.
Due to the constraints, we have that $\ve d(\ve x)\in \feas - \ve x$ and thereby 
$\ve d(\ve x) +\ve x - \ve x\ita{t} \in \feas - \ve x\ita{t}$.
Let 
$$
    (0,1] \ni 
    \sigma \ita{t} := \begin{cases}
        \min
        \left\{ 
            1, \dfrac{1}{\norm{\ve d(\ve x) + \ve x - \ve x\ita{t}}}
        \right\} 
        &\text{if $\ve d(\ve x) \ne \ve x\ita{t} - \ve x$},\\
        1&\text{else.}
    \end{cases}
$$
Then $\sigma\ita{t} \left( \ve d(\ve x) + \ve x - \ve x\ita{t}\right)$ is feasible for 
\eqref{eqn:descent_direction_problem1} at $\ve x\ita{t}$:
\begin{itemize} 
    \item $\sigma\ita{t} \left( \ve d(\ve x) + \ve x - \ve x\ita{t}\right) \in \feas - \xitat$ 
    because $\feas -\ve x\ita{t}$ is convex and 
    $\ve 0,\left( \ve d(\ve x) + \ve x - \ve x\ita{t}\right) \in \feas - \ve x\ita{t}$ as well as
    $\sigma\ita{t}\in (0,1]$.
    \item $\norm{\sigma\ita{t} \left( \ve d(\ve x) + \ve x - \ve x\ita{t}\right)}\le 1$ by the definition of $\sigma\ita{t}$.
\end{itemize}
By the definition of \eqref{eqn:descent_direction_problem1} it follows that
\begin{align}
    \max_\ell \langle \gradf_\ell (\ve x\ita{t}), \ve d(\ve x\ita t ) \rangle
        &\le 
            \sigma\ita{t} \max_\ell \langle \gradf_\ell (\ve x\ita{t}), \ve d(\ve x) + \ve x - \ve x\ita{t} \rangle 
        \nonumber
        \shortintertext{and by \cref{appendix:max_prop1}}            
    \max_\ell \langle \gradf_\ell (\ve x\ita{t}), \ve d(\ve x\ita{t}) \rangle
        &\le 
            \sigma\ita{t} \max_\ell \langle \gradf_\ell (\ve x\ita{t}), \ve d(\ve x) \rangle 
            +
            \sigma\ita{t} \max_\ell \langle \gradf_\ell (\ve x\ita{t}), \ve x - \ve x\ita{t} \rangle.
        \label{eqn:appendix1}            
\end{align}
We make the following observations:
\begin{itemize}
    \item 
    Because of 
        $\norm{ \ve d(\ve x) + \ve x - \ve x \ita{t} } \xrightarrow{t\to \infty} \norm{\ve d(\ve x)} \le 1$, 
        it follows that $\sigma \ita{t} \xrightarrow{t\to \infty} 1$.
    \item 
        Because all objective gradients are continuous, it holds for all $\ell \in \{1,\ldots, k\}$ that 
        $\gradf_\ell(\ve x\ita{t})\to \gradf_\ell(\ve x)$ and 
        because $\ve u \mapsto \max_\ell u_\ell$ is continuous as well, it then follows that 
        $$
        \max_\ell \langle \gradf_\ell (\ve x\ita{t}), \ve d(\ve x) \rangle 
        \to 
        \max_\ell \langle \gradf_\ell (\ve x), \ve d(\ve x) \rangle \quad \text{for $t\to \infty$.}
        $$
    \item
        The last term on the RHS of \eqref{eqn:appendix1} vanishes for $t\to \infty$.
\end{itemize}
By taking the limit superior on \eqref{eqn:appendix1}, we then find that 
\begin{equation}
    \limsup_{t \to \infty}
        \theta ( \ve x\ita t )
    =
    \limsup_{t \to \infty}
        \max_\ell \langle  
            \gradf_\ell (\ve x\ita{t}), \ve d(\ve x\ita{t})
        \rangle
    \le 
        \max_\ell \langle 
            \gradf_\ell (\ve x), \ve d(\ve x) 
        \rangle 
    = 
        \theta ( \ve x )
    \label{eqn:appendix2}
\end{equation}
Vice versa, we know that because of $\ve d(\ve x\ita{t})\in \feas - \ve x\ita{t}$, it holds that 
$\ve d(\ve x\ita{t}) + \ve x\ita{t} - \ve x\in \feas - \ve x$ and as above we find that 
\begin{equation} 
    \max_\ell \langle \gradf_\ell( \ve x ), \ve d(\ve x ) \rangle 
        \le 
            \lambda\ita{t}\max_\ell \langle \gradf_\ell(\ve x), \ve d( \ve x\ita{t}) \rangle 
            +
            \lambda\ita{t}\max_\ell \langle \gradf_\ell (\ve x), \ve x\ita t - \ve x\rangle
    \label{eqn:appendix6}
\end{equation}
with 
$$ 
    \lambda \ita{t} := \begin{cases}
        \min
        \left\{ 
            1, \dfrac{1}{\norm{\ve d(\ve x) + \ve x\ita{t} - \ve x }}
        \right\} 
        &\text{if $\ve d(\ve x) \ne \ve x\ita{t} - \ve x$},\\
        1&\text{else.}
    \end{cases}
$$
Again, the last term of \eqref{eqn:appendix6} vanishes in the limit so that by using the 
properties of the maximum function and the continuity of $\gradf_\ell$, 
as well as $\lambda\ita{t}\xrightarrow{t\to \infty}1$, in taking the limit inferior 
on \eqref{eqn:appendix6} we find that
\begin{align} 
    & \theta(\ve x) = \max_\ell \langle \gradf_\ell( \ve x ), \ve d(\ve x ) \rangle 
    %\hspace{-3cm}&
    %\nonumber\\
    %&
        \le 
            \liminf_{t\to \infty}
                \max_\ell \langle \gradf_\ell(\ve x), \ve d( \ve x\ita{t}) \rangle
        \nonumber\\
        &\le 
            \liminf_{t\to \infty} 
            \left[
                \left( 
                    \max_\ell \langle \gradf_\ell(\ve x), \ve d( \ve x\ita{t}) \rangle 
                    - 
                    \max_\ell \langle \gradf_\ell(\ve x\ita{t}), \ve d( \ve x\ita{t}) \rangle 
                \right)
                %\right. \nonumber\\
                %&\qquad 
                %\left. 
                +
                    \max_\ell \langle \gradf_\ell(\ve x\ita{t}), \ve d( \ve x\ita{t}) \rangle 
            \right]
        \nonumber\\
        &\le 
        \liminf_{t\to \infty} 
            \left[
                \norm{
                    \jacobianf(\ve x) 
                    -
                    \jacobianf(\ve x\ita{t})
                }
                \norm{ \ve d(\ve x\ita{t}) }
                + \max_\ell \langle \gradf_\ell(\ve x\ita{t}), \ve d( \ve x\ita{t}) \rangle 
            \right]
        \nonumber\\
        &\le 
        \liminf_{t\to \infty} 
            \max_\ell \langle \gradf_\ell(\ve x\ita{t}), \ve d( \ve x\ita{t}) \rangle 
        = 
        \liminf_{t\to \infty}
            \theta(\ve x\ita t).
        \label{eqn:appendix3}
\end{align}
Combining \eqref{eqn:appendix2} and \eqref{eqn:appendix3} shows that $\theta(\ve x\ita{t})\xrightarrow{t\to \infty} \theta(\ve x)$.
\end{proof}

\cref{theorem:omega_uniformly_continuous} claims that $\omga{\ve x}$ is uniformly continuous,
provided the objective gradients are Lipschitz. 
The implied Cauchy continuity is an important property in the convergence proof of the algorithm.
\begin{proof}[Proof of \cref{theorem:omega_uniformly_continuous}]
    We will consider the constrained case only, when $\feas$ is convex and compact 
    and show uniform continuity \emph{a fortiori} by proving that $\omga{\bullet}$ 
    is Lipschitz.
    Let the objective gradients be Lipschitz continuous.
    Then $\jacobianf$ is Lipschitz as well with constant $L > 0$. 
    Let $\ve x,\ve y\in \feas$ with $\ve x\ne \ve y$ (the other case is trivial) and 
    let again $\ve d(\ve x), \ve d(\ve y)$ be the respective optimizers.
    
    Suppose w.l.o.g. that 
    \begin{align*}
    \abs{ 
        \max_\ell 
        \langle 
            \gradf_\ell (\ve x), \ve d(\ve x)
        \rangle 
        - 
        \max_\ell 
        \langle 
            \gradf_\ell (\ve y), \ve d(\ve y)
        \rangle
    } 
    &=
        \max_\ell 
            \langle 
                \gradf_\ell (\ve x), \ve d(\ve x)
            \rangle 
        - 
        \max_\ell 
            \langle 
                \gradf_\ell (\ve y), \ve d(\ve y)
            \rangle
    \end{align*}
    If we define 
    $$
        (0,1]\ni \sigma :=
        \begin{cases} 
            \min\left\{ 1, \frac{1}{\norm{\ve d(\ve y) + \ve y - \ve x}} \right\} 
                &\text{if $\ve d(\ve y) \ne \ve x - \ve y$,}\\
            1 &\text{else,}
        \end{cases}
    $$
    then again $\sigma \left( \ve d(\ve y) + \ve y - \ve x\right)$ 
    is feasible for \eqref{eqn:descent_direction_problem1} 
    at $\ve y$.
    Thus, 
    \begin{align}
        %\abs{ 
        \max_\ell 
        \langle 
            \gradf_\ell (\ve x), \ve d(\ve x)
        \rangle 
        - 
        \max_\ell 
        \langle 
            \gradf_\ell (\ve y), \ve d(\ve y)
        \rangle
    %}
    \hspace{-3cm}&
    \nonumber\\
    &\stackrel{\text{df.}}\le 
    \max_\ell 
        \langle 
            \gradf_\ell (\ve x), \sigma \left( \ve d(\ve y) + \ve y - \ve x\right)
        \rangle 
        - 
        \max_\ell 
            \langle 
                \gradf_\ell (\ve y), \ve d(\ve y)
            \rangle
    \nonumber\\
    &\le 
        \norm{
            \sigma \jacobianf(\ve x) \left( \ve d(\ve y) + \ve y - \ve x\right)
            -
            \jacobianf(\ve y) \ve d(\ve y)
        }
    \nonumber\\
    &\stackrel{\sigma \le 1}\le 
        \norm{ 
            \sigma \jacobianf(\ve x) - \jacobianf(\ve y)
        }
        \norm{ \ve d( \ve y )}
        + 
        \norm{\jacobianf(\ve x)} \norm{\ve x - \ve y},
    \label{eqn:appendix4}
    \end{align}
    where we have again used \cref{appendix:max_prop2} for the second inequality.
    We now investigate the first term on the RHS.
    Using $\norm{\ve d(\ve y)} \le 1$ and adding a zero, we find 
    \begin{align} 
    \norm{ 
        \sigma \jacobianf(\ve x) - \jacobianf(\ve y)
    }
    \norm{ \ve d( \ve y ) }
    % &\le 
    %     \norm{ 
    %         (1-\sigma)\jacobianf(\ve x) + \sigma \jacobianf(\ve x) - \jacobianf(\ve y) - (1-\sigma) \jacobianf(\ve x)
    %     }
    % \nonumber
    % \\
    &\le 
    \norm{ 
        \jacobianf(\ve x) - \jacobianf(\ve y) - (1-\sigma) \jacobianf(\ve x)
    }
    \nonumber
    \\
    &
    \le 
        L \norm{ 
            \ve x - \ve y
        }
        + (1-\sigma) \norm{\jacobianf(\ve x)}.
    \label{eqn:appendix5}
    \end{align}
    Furthermore, $\norm{\ve d(\ve y) + \ve y - \ve x} \le 1 + \norm{\ve y - \ve x}$ implies 
    $1/( 1 + \norm{\ve y - \ve x} ) \le \sigma$ and 
    $$
        1 - \sigma \le 1 - \frac{1}{ 1 + \norm{\ve y - \ve x} } = 
        \frac{ \norm{\ve y - \ve x} }{ 1 + \norm{\ve y - \ve x} }
        \le \norm{\ve y - \ve x}.
    $$
    We use this inequality and plug \eqref{eqn:appendix5} into \eqref{eqn:appendix4} to obtain 
    \begin{align*}
        %\abs{ 
        \max_\ell 
        \langle 
            \gradf_\ell (\ve x), \ve d(\ve x)
        \rangle 
        - 
        \max_\ell 
        \langle 
            \gradf_\ell (\ve y), \ve d(\ve y)
        \rangle
    %}  
        &\le
            L\norm{ \ve x - \ve y} 
            + 
            2\norm{\jacobianf(\ve x)} \norm{ \ve x - \ve y} 
        \\
        &\le (L + 2D) \norm{ \ve x - \ve y},
    \end{align*}
    with $D = \max_{\ve x \in \feas} \norm{ \jacobianf(\ve x)}$ which is well defined 
    because $\feas$ is compact and $\norm{\jacobianf({\bullet})}$ is continuous.
\end{proof}

%\subsection{Sufficient Decrease}
%\label{appendix:section_sufficient_decrease}

\subsection{Modified Criticality Measures}
\label{appendix:modomega}

\begin{proof}[Proof of \cref{theorem:omegaconst_mod}]
    There are two cases to consider:
    \begin{itemize}
      \item If $\omegamt{ \xitat } \ge \omga{ \xitat }$ then
      $$
        \left|
            \omegamt{ \xitat } - \omga{  \xitat  }
        \right|
        =
        \omegamt{ \xitat } - \omga{  \xitat  }
        \le
        \omegaconstant \omegamt{  \xitat  }.
      $$
      Now
      \begin{equation*}
        {
        \left|
            \modomegamt{ \xitat } - \modomega{  \xitat  }
        \right|
        \in }
        \left\{
          \begin{aligned}
            &\phantom{{}\le{}} \omegamt{ \xitat } - \omga{  \xitat  }
            \\
            1 - \omga{  \xitat  }  &\le \omegamt{ \xitat } - \omga{  \xitat  }\\
            1 - 1 &= 0
          \end{aligned}
        \right\}
        \le \omegaconstant \omegamt{  \xitat  }.
      \end{equation*}
      \item The case $\omga{ \xitat } < \omegamt{ \xitat }$ can be shown similarly.
    \end{itemize}
\end{proof}

\begin{proof}[Proof of \cref{theorem:omegaconst}]
    Use \cref{theorem:omegaconst_mod} and then investigate the two possible cases:
    \begin{itemize}
      \item If $\modomegamt{ \xitat  }\ge \modomega{ \xitat }$, then the first inequality follows because of $1 \ge 1/(1 + \omegaconstant)$.
      \item If $\modomegamt{ \xitat } < \modomega{ \xitat }$, then
        $\modomega{ \xitat } - \modomegamt{ \xitat } \le \omegaconstant \modomegamt{  \xitat },$
      and again the first inequality follows.
    \end{itemize}
\end{proof}

\section{Pascoletti-Serafini Step}
\label[Appendix]{appendix:section_steps}
%As in the scalar case (see \cite{conn:trm_framework}) we do not have to use the steepest-descent direction
%but can employ other steps as long as they provide a bound like
%\eqref{eqn:sufficient_decrease_prototype_modomega}.
%\spfoot{Die Gl.\ kommt ers tviel später. Lässt sich dieser Hinweis vermeiden (bzw.\ einfach streichen)?}.
%That is, we can relate them to the criticality measure induced by the surrogate steepest-descent direction.

%\subsubsection{Pascoletti-Serafini Steps}
One example of an alternative descent step $\stept\in \mathbb{R}^n$ is given in \cite{thomann}.
\citet*{thomann} leverage the Pascoletti-Serafini scalarization to define local 
subproblems that guide the iterates towards the (local) model ideal point.
To be precise, it is shown that the trial point $\xtrialt$ can be computed as the solution to
\begin{equation}
  \begin{aligned}
    \min_{\tau \in \mathbb{R}, \ve x \in \ballt} \tau \quad \text{ s.t. }
    \vemt(\xitat) + \tau\ve r\itat - \vemt(\ve x) \ge \ve 0,
  \end{aligned}
\label{eqn:pascoletti_serafini_m}
\end{equation}
where $\ve r\itat =\vemt(\xitat) - \ve i_{\mathrm m}\ita t\in \mathbb{R}^k_{\ge 0}$ is the direction vector pointing from the local model ideal point
\begin{equation}
  \ve i_{\mathrm{m}}\ita t
  =
  \begin{bmatrix}
    i_1\itat,
    \ldots,
    i_k\itat
  \end{bmatrix}^T,
  \,\text{ with
    $i_\ell\itat = \min_{\ve x\in \feas} \mitat_\ell(\ve x)$ for $\ell = 1, \ldots, k,$}
\label{eqn:local_model_ideal_point}
\end{equation}
to the current iterate value.
If the surrogates are linear or quadratic polynomials and the trust region use a $p$-norm 
with $p\in \{1,2,\infty\}$ these sub-problems are linear or quadratic programs.

A convergence proof for the unconstrained case is given in \cite{thomann}.
It relies on a sufficient decrease bound similar to \eqref{eqn:sufficient_decrease_prototype_modomega}.
However, it is not shown that $\sufficientconstant \in (0,1)$ exists 
independent of the iteration index $t$ but stated as an assumption.\\
Furthermore, constraints (in particular box constraints) are integrated into 
the definition of $\omga{\bullet}$ and $\omegamt{\bullet}$ using an active set strategy 
(see \cite{thomann_diss}).
%instead of the projection approach used in this article.
Consequently, both values are no longer Cauchy continuous.
We can remedy both drawbacks by relating the (possibly constrained) Pascoletti-Serafini 
trial point to the strict modified Pareto-Cauchy point in our projection framework.
To this end, we allow in \eqref{eqn:pascoletti_serafini_m} 
and \eqref{eqn:local_model_ideal_point} any feasible set fulfilling 
\cref{assumption:feasible_set_compact_convex}.
Moreover we recite the following assumption:
\begin{Assumption}[Assumption 4.10 in \cite{thomann}]
  There is a constant $\rconstant\in(0,1]$ so that if $\xitat$ is not Pareto-critical, 
  the components $r_1\itat, \ldots, r_k\itat,$ of $\ve r\itat$ satisfy
  %\begin{equation*}
    $
    \dfrac{
      \min_\ell r_\ell\itat
    }{
      \max_\ell r_\ell\itat
    } \ge \rconstant.$
    %\end{equation*}
  \label{assumption:rconstant}
\end{Assumption}
The assumption can be justified because $r_\ell\itat > 0$ if $\xitat$ is not critical 
and $r_\ell\itat$ can be bounded above and below by expressions involving 
$\omegamt{\bullet}$, see \cref{remark:minimizer} and \cite[Lemma 4.9]{thomann}.
We can then derive the following lemma:
\begin{Lemma}
  Suppose \cref{assumption:feasible_set_compact_convex,assumption:f_m_twice_continuously_differentiable,assumption:rconstant} hold.
  Let $(\tau^+,\xtrialt)$ be the solution to \eqref{eqn:pascoletti_serafini_m}.
   Then there exists a constant $\tildesufficientconstantm\in (0,1)$ such that it holds
  $$
  \Phimt(\xitat) - \Phimt(\xtrialt)
  \ge
  \tildesufficientconstantm
  \omegamt{  \xitat  } \min
  \left\{
    \frac{\omegamt{ \xitat }}{\normconst\hessboundmt}, \radiust, 1
  \right\}.
  $$
  \label{lemma:ps_sufficient_decrease}
\end{Lemma}

\begin{proof}
If $\xitat$ is critical for \eqref{eqn:mopmt}, then $\tau^+=0$ and $\xtrialt = \xitat$ and 
the bound is trivial \cite{eichfelder:adaptive}.
Otherwise, we can use the same argumentation as in \cite[Lemma 4.13]{thomann} to show that 
for the strict modified Pareto-Cauchy point $\strictx$ it holds that
\begin{align*}
\Phimt(\xitat) - \Phimt(\xtrialt)
    &\ge \rconstant \min_\ell\left\{ \mitat_\ell(\xitat) - \mitat_\ell(\strictx)\right\}
%     \shortintertext{and from \cref{theorem:strict_sufficient_decrease}}
% \Phimt(\xitat) - \Phimt(\xtrialt)
%     &\ge
%         \underbrace{
%             \rconstant
%             \sufficientconstantm
%         }
%         _{=:\tildesufficientconstantm \in (0,1)}
%             \omegamt{  \xitat  } \min
%             \left\{
%                 \frac{\omegamt{ \xitat }}{\normconst\hessboundmt}, \radiust, 1
%             \right\}
%     .
\end{align*}
and the final bound follows from \cref{theorem:strict_sufficient_decrease} with the 
new constant $\tildesufficientconstantm = \rconstant\sufficientconstantm$.
\end{proof}

% The Pascoletti-Serafini step is particularly suited if
% \begin{itemize}
% 	\item the trust regions are determined using $\norm{\bullet}_p$ with $p\in \{1,2,\infty\}$,
%  	\item $\feas = \mathbb{R}^n$ or $\feas$ is constrained by linear or quadratic constraints only and
% 	\item the surrogates $\vemt$ are multivariate polynomials of degree at most 2.
% \end{itemize}
% Then \eqref{eqn:pascoletti_serafini_m} is a linear or quadratic program and can
% be solved with reasonable numerical effort.
% Further, instead of local ideal points global utopia vectors can be provided to guide the descent.
%\input{trash/nlcg.tex}
%\input{trash/appendix_gradient_related.tex}

%%%%%%%%%%%%%%%%%%%%%%%%%%%%%%%%%%%%%%%%%%
\end{paracol}
\reftitle{References}


%\nocite{*}
\externalbibliography{yes}
\bibliography{sources.bib}

% If authors have biography, please use the format below
%\section*{Short Biography of Authors}
%\bio
%{\raisebox{-0.35cm}{\includegraphics[width=3.5cm,height=5.3cm,clip,keepaspectratio]{Definitions/author1.pdf}}}
%{\textbf{Firstname Lastname} Biography of first author}
%
%\bio
%{\raisebox{-0.35cm}{\includegraphics[width=3.5cm,height=5.3cm,clip,keepaspectratio]{Definitions/author2.jpg}}}
%{\textbf{Firstname Lastname} Biography of second author}

% The following MDPI journals use author-date citation: Arts, Econometrics, Economies, Genealogy, Humanities, IJFS, JRFM, Laws, Religions, Risks, Social Sciences. For those journals, please follow the formatting guidelines on http://www.mdpi.com/authors/references
% To cite two works by the same author: \citeauthor{ref-journal-1a} (\citeyear{ref-journal-1a}, \citeyear{ref-journal-1b}). This produces: Whittaker (1967, 1975)
% To cite two works by the same author with specific pages: \citeauthor{ref-journal-3a} (\citeyear{ref-journal-3a}, p. 328; \citeyear{ref-journal-3b}, p.475). This produces: Wong (1999, p. 328; 2000, p. 475)

%%%%%%%%%%%%%%%%%%%%%%%%%%%%%%%%%%%%%%%%%%
%% for journal Sci
%\reviewreports{\\
%Reviewer 1 comments and authors’ response\\
%Reviewer 2 comments and authors’ response\\
%Reviewer 3 comments and authors’ response
%}
%%%%%%%%%%%%%%%%%%%%%%%%%%%%%%%%%%%%%%%%%%
\end{document}